\documentclass{article}

\usepackage{PRIMEarxiv}
\usepackage[utf8]{inputenc}
\usepackage[utf8]{inputenc}
\usepackage{textcomp}
\usepackage{multirow}
\usepackage{graphicx}
\usepackage{amsmath}
\usepackage[version=4]{mhchem}
\usepackage{siunitx}
\usepackage{longtable,tabularx}
\usepackage{mathalpha}
\usepackage{bm}
\usepackage{adjustbox}
\usepackage{footnpag}			  
\usepackage{acronym}
\usepackage{caption}
\usepackage{subcaption}
\usepackage{siunitx}
\usepackage{booktabs}%
\usepackage{hyperref}
\usepackage{verbatim}
\usepackage{algorithm}
\usepackage{algpseudocode}
\usepackage{mathtools}
\usepackage[thinc]{esdiff}
\usepackage{capt-of}
\usepackage{longtable,tabularx}
\usepackage{caption}
\usepackage{subcaption}
\usepackage{arydshln}
\usepackage{setspace}
\usepackage{enumitem}
\usepackage{tikz}
\usetikzlibrary{shapes, arrows}
 \usetikzlibrary{calc} 
 \usetikzlibrary{fit}
 \usepackage[nocomma]{optidef}
 \usetikzlibrary{automata,positioning}
\usepackage{tablefootnote}
 \usepackage[nocomma]{optidef}
 \usetikzlibrary{automata,positioning}
\usepackage{tablefootnote}

\newcommand{\mat}{\begin{bmatrix}}
\newcommand{\matf}{\end{bmatrix}}
\newcounter{equationset}
\usepackage{mathrsfs}
\usepackage{amssymb}
\title{Trajectory Design and Guidance for Far-range Proximity Operations of Active Debris Removal Missions with Angles-only Navigation and Safety Considerations
}

\author{
  Minduli C. Wijayatunga \\
  Phd Candidate\\
  Te Pūnaha Ātea - Space Institute \\
  The University of Auckland \\
  Auckland\\ \\
   \And
  Roberto Armellin\\
  Professor \\
  Te Pūnaha Ātea - Space Institute \\
  The University of Auckland \\
  Auckland\\ \\
       \And
Harry Holt \\
  Research Fellow \\
  Te Pūnaha Ātea - Space Institute \\
  The University of Auckland \\
  Auckland\\\\}

\begin{document}
\maketitle

\begin{abstract}
Observability of the target, safety, and robustness are often recognized as critical factors in ensuring successful far-range proximity operations. The application of angles-only (AO) navigation for proximity operations is often met with hesitancy due to its inherent limitations in determining range, leading to issues in target observability and consequently, mission safety. However, this form of navigation remains highly appealing due to its low cost. This work employs Particle Swarm Optimization (PSO) and Reinforcement Learning (RL) for the design and guidance of such far-range trajectories, assuring observability, safety and robustness under angles-only navigation. Firstly, PSO is used to design a nominal trajectory that is observable, robust and safe. Subsequently, Proximal Policy Optimization (PPO), a cutting-edge RL algorithm, is utilized to develop a guidance controller capable of maintaining observability while steering the spacecraft from an initial perturbed state to a target state. The fidelity of the guidance controller is then tested in a Monte-Carlo (MC) manner by varying the initial relative spacecraft state.  The observability of the nominal trajectory and the perturbed trajectories with guidance are validated using an Extended Kalman Filter (EKF). The perturbed trajectories are also shown to adhere to the safety requirements satisfied by the nominal trajectory.
Results demonstrate that the trained controller successfully determines maneuvers that maintain observability and safety and reaches the target end state. 

\end{abstract}

\section*{Nomenclature}


{\renewcommand\arraystretch{1.0}
\noindent\begin{longtable*}{@{}l @{\quad=\quad} l@{}}
$PSO$  & Particle Swarm Optimization \\
$EKF$ &   Extended Kalman Filter \\
$PPO$& Proximal Policy Optimization \\
$RL$ & Reinforcement Learning \\
$MC$ & Monte Carlo \\
$g_0$ & Gravitational acceleration (9.80665 m/s$^2$) \\
$\mu$& Gravitational parameter ($\SI{3.986e5}{\kilo\meter\cubed\per\second\squared}$)\\ 
\multicolumn{2}{@{}l}{Subscripts}\\
$bal$ & Absolute ballistic trajectory of the servicer \\
$i$ & Segment index\\
$mp$	& Midpoint of a segment \\
\multicolumn{2}{@{}l}{Superscripts}\\
$d$ & Debris orbit parameter\\
$rel$ & Relative orbit parameter\\
\end{longtable*}}

\section{Introduction}
\label{intro}
RPO missions generally refer to any mission that performs orbital maneuvers in which two spacecraft arrive in the same orbit and approach at a close distance \cite{YUAN2022107812}. This rendezvous may or may not be followed by a docking procedure. Space missions with RPO are becoming increasingly commonplace due to enhanced interest in space target removal \cite{ADRMW, convexMW}, on-orbit servicing \cite{CHAI2022107527,RODRIGUES2022107865}, and asteroid inspection/mining/deflection \cite{HU2024109316,CHU2024108877} missions. The far-range phase of an RPO mission starts from the point where relative navigation measurements can be made and ends at a point at least \SI{1}{\kilo\meter} away from the target \cite{BORELLI202333}. Safety, observability, and robustness are critical drivers for these operations, especially when targetting uncoorperative objects \cite{BORELLI202333}. \par 

 Although safety in far-range approaches to cooperative targets has a strong heritage that stems from the Automated Transfer Vehicle (ATV) missions \cite{ganet2002atv}, the Orbital Express and Engineering Test Satellite No. 7 (ETS-VII) demonstration missions \cite{mulder2008orbital,kawano2001result}, safe approaches to uncooperative target are less explored but require a higher level of safety \cite{gaias2014angles}. The use of centered and off-centered Walking Safety Ellipse (WSE) trajectories for spacecraft safety was explored by D. Gaylor in \cite{6dc118e90c5e46b1938eaf2e45ecf68a}. WSE can be applied for uncooperative spacecraft RPO as it ensures that the trajectory never crosses the velocity vector of the target, thereby assuring passive (PAS) safety. However, WSE trajectories can be fuel-expensive \cite{MinduliThesis}. In \cite{Barbee2011AGA}, co-elliptic trajectories have been explored for RPO, where the chaser is situated in the same orbital plane as the target while having a different semi-major axis. However, this does not assure PAS safety when the chaser approaches the target's orbit. Following the development of Relative Orbital Elements (ROEs) in \cite{simonethesis}, spiral approaches with Eccentricity/Inclination (E/I) separation for ensuring PAS safety were developed  \cite{d2006proximity}. E/I separation works by ensuring that the chaser and the target maintain a minimum separation distance in the radial-cross-track plane, as the state uncertainties in those directions are significantly less compared to the along-track direction \cite{d2006proximity}. This method has been demonstrated in orbit via the Autonomous Vision Approach Navigation and Target Identification (AVANTI) experiment \cite{gaias2016design} and the Restore-L servicing mission \cite{reed2016restore} on uncooperative and semi-cooperative targets. This work uses a minimum radial-normal distance metric derived from the E/I separation formulation to enforce mission safety. 

 Accurate navigation is another critical requirement of successful RPO missions. Although {AON} is advantageous in terms of cost and simplicity, determining the entire state of the target only based on angular measurements can result in difficulties in determining the range \cite{obsangle}. Woffinden \cite{woffinden2009observability} first presented the criteria for {AON} observability using a geometric interpretation, confirming that altering the natural line of sight measurement profile could guarantee observability. He also proposed several detectability metrics that can characterize the degree of observability in {AON} measurements under the Hill-Clohessy Whitsire (HCW) equations \cite{doi:10.2514/1.45006}, optimizing trajectories for observability. Gaias et al. proposed yet another metric of observability involving the maximum amount of change in the measurement angles between forced and natural motion trajectories \cite{gaias2014angles}. Franquiz et al. also adopted a similar metric in \cite{obsangle} to develop sequences of two-burn impulsive maneuvers that improve observability. This work adopts a similar concept: the dot product between the natural vs. forced measurement profile is used to estimate observability.

Robustness and accurate guidance during far-range operations are critical for the success of RPO missions. Luo et al. introduced a closed-loop guidance scheme for AON using HCW equations \cite{LUO201691}. However, this approach does not account for safety considerations. Hablani et al. proposed a guidance scheme for far-range trajectories based on the glidescope method \cite{doi:10.2514/2.4916}, but it similarly overlooks navigation and safety criteria. In another study, Yang et al. developed a guidance framework for autonomous collision avoidance under uncertainties in rendezvous missions \cite{YANG2024109547}. Li et al. presented a nonlinear model predictive control guidance method for spacecraft rendezvous using AON, formulating it as a closed-loop optimal control problem \cite{LI2017236}. Guo et al. addressed the challenge of autonomous docking with uncooperative targets, introducing a guidance scheme that incorporates safety constraints \cite{GUO2021106380}. Additionally, Yuan et al. developed a reinforcement learning-based guidance scheme for far-range approaches, integrating AON navigation and approach constraints while also utilizing HCW equations \cite{YUAN2022107812}.

In contrast, this work utilizes Particle Swarm Optimization (PSO) and Reinforcement Learning ({RL}) to develop far-range trajectories taking all the discussed criteria into account. Evolutionary algorithms such as {PSO} were first developed to address the limitations of gradient-based methods, including their local nature and requirement for a suitable initial guess in the region of convergence \cite{pontani2012particle}.   
 {PSO} is an evolutionary algorithm first developed by Kennedy and Eberhart in 1995 \cite{eberhart1995new,kennedy1995particle}. Similar to all evolutionary algorithms, {PSO} is also motivated by nature by utilizing a population of individuals that symbolize potential solutions to the target problem. The best solution is then obtained via collaboration and competition in the population. 
In particular, it mimics the motion of foraging in a flock of birds and their information sharing mechanism that affects the behavior of the whole flock \cite{pontani2012particle}.  
The use of {PSO} in RPO missions is commonly seen in the literature, starting with the use of {PSO} for solving time-limited orbit transfers and impulsive interplanetary transfers in 2006 by Bessete and Spencer \cite{bessette2006optimal, bessette2006optimal}. Zhu et al. utilized {PSO} to optimize low thrust trajectory for asteroid exploration \cite{zhu2009trajectory}. Pontani et al. utilized {PSO} for the optimization of impulsive and low thrust trajectories using nonlinear equations and {HCW} equations \cite{pontani2012particle,doi:10.2514/1.A32402}. However, using {PSO} for nominal far-range trajectory design under safety and observability constraints remains unexplored to the authors' knowledge. 

The origins of {RL} lies in the development of optimal control and Markov decision processes in the 1960s \cite{BELLMAN1958228}, as well as in the learning by trial and error work done on the psychology of animal training   \cite{montague1999reinforcement}. 
Neural networks have become popular for policy representation in {RL} due to their generalization capabilities and scalability \cite{HORNIK1989359}. {RL} has frequently been applied to space-related tasks in the last decade. Chan et al. used it for spacecraft map generation while orbiting small bodies \cite{chan2019autonomous}, Willis and Izzo utilized {RL} to develop spacecraft orbit control laws in unknown gravitational fields \cite{willis2016reinforcement} and Zavoli and Federici has used it for the trajectory design of interplanatary missions \cite{doi:10.2514/1.G005794}. {RL} has also been used for spacecraft orbital transfers \cite{lafarge2020guidance} and to perform guidance and control for pinpoint planetary landing \cite{scorsoglio2020image,gaudet2014adaptive}. Building off this work, in \cite{doi:10.2514/1.A34838}, Ulrich et al. developed an {RL}-based guidance scheme for spacecraft pose tracking and docking.  The work presented in \cite{FEDERICI2022129} developed an adaptive guidance scheme for conducting multiple space rendezvous missions, and illustrateed that it can reconstruct optimal control solutions in the analyzed cases. Similarly, \cite{doi:10.2514/1.G00683} presented the successful use of RL for guidance during the final phase of an asteroid impactor mission.  RL is noted to be advantageous for spacecraft guidance as it allows users to develop guidance laws for complex tasks with little to no user input \cite{doi:10.2514/1.A34838}.

This work develops a robust guidance scheme that meets safety and observability criteria by training a controller via {RL} in a stochastic environment.  The overall methodology involved in this work can be divided into two stages. 
\begin{itemize}
    \item \textbf{Stage 1:} {PSO}-based nominal trajectory design considering safety and observability.
    \item \textbf{Stage 2:} {RL}-based guidance scheme development considering initial state errors, thrust errors, observability, and safety. 
\end{itemize}

The novelty of this work is the development of a two-stage, cohesive tool for trajectory design and guidance for far-range approaches, considering robustness, safety, and observability. Firstly, {PSO} is used to design a nominal trajectory that is observable and safe. Subsequently, Proximal Policy Optimization ({PPO}), a cutting-edge {RL} algorithm, is utilized to develop a guidance scheme capable of maintaining safety and observability while steering the spacecraft from an initial perturbed state to a target state. The problem formulation is another novelty of this work, where, as both PSO and RL tend to struggle with constraint satisfaction \cite{Paduraru2021ChallengesOR,Gad2022ParticleSO}, the far-range problem is formulated so that the final state condition is always satisfied.  The fidelity of the guidance scheme is then tested in a Monte-Carlo ({MC}) manner by varying the initial relative spacecraft state and introducing thrust errors. The observability of the nominal and perturbed trajectories with guidance is validated using an Extended Kalman Filter ({EKF}). The PAS safety and point-wise (PWS) safety of the generated trajectories is validated by checking the minimum distance to the target. 

The remainder of this paper is organized as follows. Section 2 delves into the proposed methodology, which is divided into two parts: PSO-based nominal trajectory design and RL-based guidance design. Section 3 presents the results of the study, where the developed methodology is applied to a 4 h far-range test case previously presented in literature. It includes an analysis of the nominal trajectories and the performance of the guidance controller in terms of fuel consumption, safety and observability. Finally, Section 4 offers concluding remarks.

\section{Methodology}
As mentioned in the introduction, the methodology of far-range trajectory design and guidance developed in this paper is divided into two phases. The first phase involves the generation of a {PSO}-based nominal trajectory. The second phase entails the design of the {RL}-based guidance scheme.

\subsection{Particle Swarm Optimization-based Nominal Trajectory Design } \label{propnode}
The {PSO} formulation for the design of the nominal trajectory optimizes the nominal $\Delta v$ impulses (denoted $\Delta \overline{\pmb{v}}$ ) of the chaser trajectory so that the following objective function can be minimized. 
\begin{equation}\label{Gmain}
   \min_{ \Delta  \overline{\pmb{v}} } G =   \gamma_1 G_{ \Delta  \overline{\pmb{v}}}  + \gamma_2 G_{\text{obs}} + \gamma_3 G_{\text{safety}} 
\end{equation}


subject to
\begin{equation}
-  \Delta \overline{{v}}_{\text{lim}}  \leq  \Delta \overline{{v}}_k \leq  \Delta \overline{{v}}_{\text{lim}}  \ \text{for} \ k = [1,...,3n]
\end{equation}

$n$ is the total number of nodes where impulses are applied, equally spaced in time. The path between any two successive nominal impulses is called a segment. $G$ indicates the objective function, and $\gamma$ indicates the weights associated with different objective function components.  $ G_{ \Delta  \overline{\pmb{v}}},G_{\text{obs}}$ and $G_{\text{safety}} $  denote the contribution to the objective function from the $\Delta v$ consumption, observability and safety, respectively.  The weights $\gamma$ are set by optimizing each component of the cost function individually and letting $\gamma$ be the inverse of the optimal individual cost. $\Delta \overline{{v}}_{\text{lim}}$ is the user-defined limit for each of the $\Delta \overline{\pmb{v}}$s.  

The chaser and the target must be propagated in time during the far-range operations. Due to the short {TOF} and the relative motion involved, only Keplerian dynamics are considered for the far-range phase.  Since the chaser and the target are at similar altitudes during far-range operations, they experience perturbations of comparable magnitudes. Consequently, perturbative effects such as atmospheric drag and $J_2$ can be assumed to have a negligible impact on relative motion. To account for any residual impacts these perturbations may have on the relative motion, initial state uncertainties are introduced later in this chapter.

The nominal state of the chaser $ \overline{\pmb{x}}(t) =   [\overline{\pmb{r}}(t); \overline{\pmb{v}}(t)]$ and the nominal state of the target $ \pmb{x}^T (t) =   [\pmb{r}^T(t); \pmb{v}^T(t)] $ follow the dynamics given by 
\begin{align}
\label{eq21}
\ddot{\overline{\pmb{x}}}(t) &=  f(\overline{\pmb{x}}(t) + \Delta  \overline{\pmb{v}})\\ 
 \ddot{{\pmb{x}}}^T(t) &= f(\pmb{x}^T(t))
\end{align}
where $ f(\pmb{x} )= -\mu_{\oplus} {\pmb{r}}/{r^3}$ and $\mu_{\oplus}$ is the gravitational parameter of Earth. 

Note that at the start time of the far-range operations $t_0$, $\overline{\pmb{x}}(t_0) = \overline{\pmb{x}}_0$ and $\overline{\pmb{x}}^T(t_0) = \overline{\pmb{x}}^T_0 $. The ECI target state for the chaser at the end time $t_f$ is denoted $\overline{\pmb{x}}(t_f) = \overline{\pmb{x}}_f$. Also, note that $\overline{\Box}$ denotes quantities related to the nominal trajectory.  Practically, as only Keplerian dynamics are considered, the time-consuming process of integrating the dynamics can be avoided by propagating in Keplerian elements and converting them back to Cartesian. As this way of propagation only requires a linear update in the mean anomaly, it is significantly faster.  The relative nominal state of the chaser with respect to the target ($\overline{\boldsymbol{x}}^{rel}$) in Radial-Normal-Tangential (RTN) coordinates is used later in the determination of the objective function components. The RTN axes are defined as 
\begin{enumerate}
    \item x-axis (R): outward along the the vector from Earth to the target. 
     \item z-axis (N): Along the angular momentum direction of the target orbit around Earth. 
    \item y-axis (T): completes the right-handed coordinate system, lies perpendicular to the x-axis in the orbital plane. 
    
\end{enumerate}
$\overline{\boldsymbol{x}}^{rel}$ can be obtained via the following conversion, which provides the relative state of the chaser in target-centric RTN coordinates. 
\begin{equation}\label{conv1}
    \overline{\boldsymbol{x}}^{rel} = [\overline{\boldsymbol{r}}^{rel} ; \overline{\boldsymbol{v}}^{rel} ] = \begin{bmatrix} \boldsymbol{C}_{ECI \rightarrow RTN} &  \boldsymbol{0}_{3 \times 3}  \\ \boldsymbol{\omega} \boldsymbol{C}_{ECI \rightarrow RTN}  &  \boldsymbol{C}_{ECI \rightarrow RTN}   \end{bmatrix} [\overline{\boldsymbol{x}} - \boldsymbol{x^T}]
\end{equation}
  where \begin{equation}
       \boldsymbol{C}_{ECI \rightarrow RTN} =  \begin{bmatrix} \frac{\boldsymbol{r}^T}{\parallel \boldsymbol{r}^T \parallel} \\  \frac{\boldsymbol{r}^T \times \boldsymbol{v}^T }{\parallel \boldsymbol{r}^T \times \boldsymbol{v}^T \parallel}  \times \frac{\boldsymbol{r}^T}{\parallel \boldsymbol{r}^T \parallel} \\ \frac{\boldsymbol{r}^T \times \boldsymbol{v}^T }{\parallel \boldsymbol{r}^T \times \boldsymbol{v}^T \parallel}   \end{bmatrix} 
  \end{equation}
  and 
 \begin{equation}
      \boldsymbol{\omega} = \frac{\parallel  \boldsymbol{r}^T \times \boldsymbol{v}^T \parallel }{{r^T}^2 }\begin{bmatrix} 0 & 1 & 0 \\ -1 & 0 & 0 \\ 0& 0 & 0 \end{bmatrix}
 \end{equation}
The starting relative state is denoted as $ \overline{\boldsymbol{x}}^{rel}_0$ and the target relative state is labeled $\overline{\boldsymbol{x}}^{rel}_f$. 

\subsubsection{Nominal $\Delta \overline{v}$ Metric}
\begin{figure}[hbt!]
      \centering
      \includegraphics[width=0.8\linewidth]{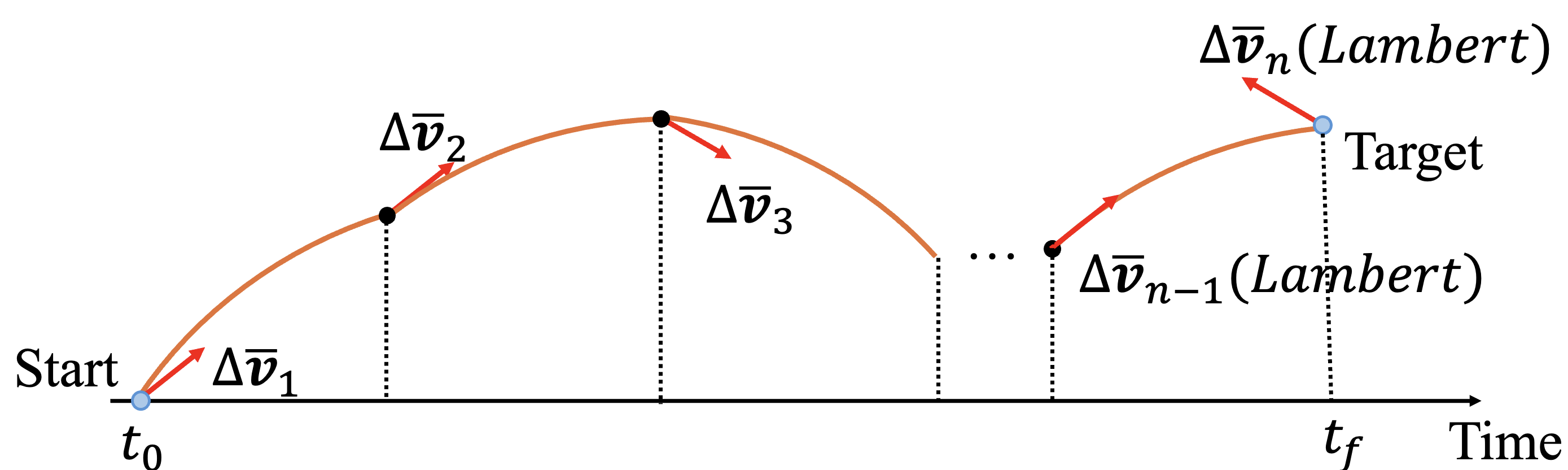}
      \caption{Nominal $\Delta \overline{v}$ metric: implementation of nominal $\Delta v$s to reach the target.}
      \label{fig:dvs}
  \end{figure}

The nominal $\Delta \overline{v}$ metric corresponds to the fuel consumption of the trajectory. 
  Note that the penultimate ($n-1^{\text{th}}$) and the ultimate ($n^{\text{th}}$) $\Delta  \overline{\pmb{v}}$ are calculated via the Lambert's method \cite{Vallado2013} such that the target ECI state $\pmb{x}_f$ is always reached during the {PSO} iterations. As PSO  and RL are known to struggle with constraint satisfaction \cite{Gad2022ParticleSO,Paduraru2021ChallengesOR}, this stops them both from encountering infeasible solutions. 
  This formulation is illustrated in Fig. \ref{fig:dvs}. 

 Hence, the objective function contribution of $\Delta  \overline{\pmb{v}}$s is then obtained by summing up all the ${\Delta \overline{\bm{v}}}$ magnitudes as follows. 
  
  \begin{equation}\label{Gdv}
      G_{\Delta  \overline{\pmb{v}}} = \sum^{n-2}_{i =0} \parallel \Delta  \overline{\pmb{v}}_i \parallel + 
     \parallel  \Delta \overline{\pmb{v}}_{n-1} (Lambert) \parallel +  \parallel  \Delta \overline{\pmb{v}}_{n} (Lambert) \parallel
  \end{equation}


  \subsubsection{Observability Metric}\label{OM}
\begin{figure}[hbt!]
    \centering
    \includegraphics[width=0.8\linewidth]{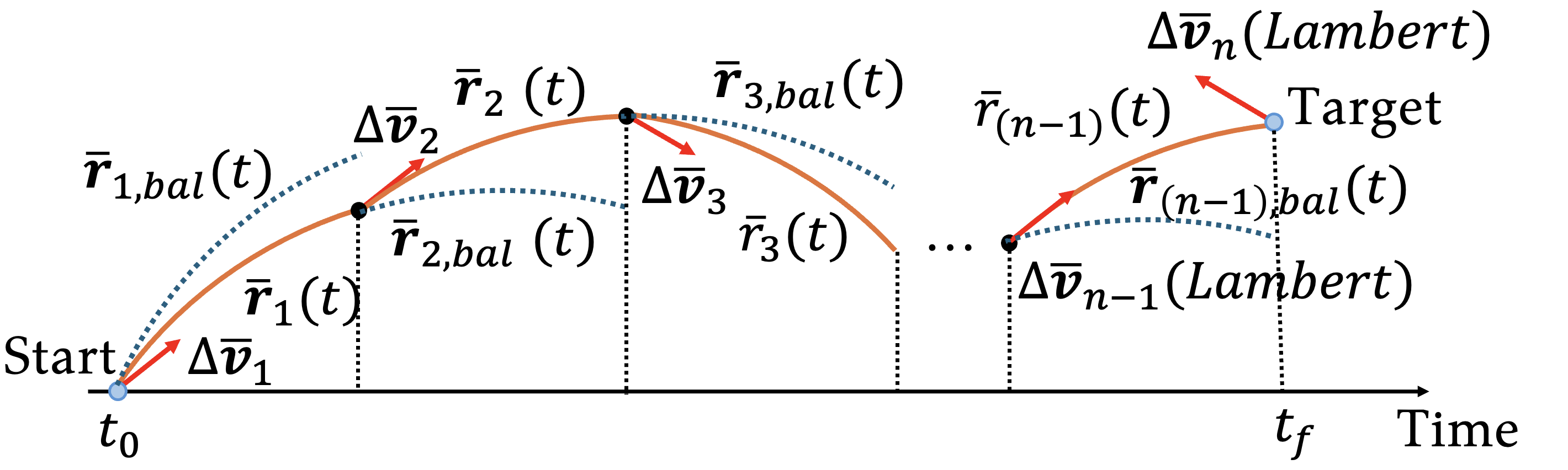}
    \caption{Observability metric: $\bm{\bar{r}}_1(t),\bm{\bar{r}}_2(t),\bm{\bar{r}}_3(t),\bm{\bar{r}}_{(n-1)}(t)$ are the maneuvred nominal trajectory position vectors, $\bm{\bar{r}}_{1,bal}(t),\bm{\bar{r}}_{2,bal}(t),\bm{\bar{r}}_{3,bal}(t),\bm{\bar{r}}_{(n-1),bal}(t)$ are the ballistic trajectory position vectors. }
    \label{obscalc}
\end{figure}

When relying only on {AON}, ensuring that a generated profile of angle measurements can fully determine the relative state of the target from the chaser spacecraft is crucial. This criterion can be achieved by maneuvering the chaser spacecraft to optimize the observability of the measurement profile obtained. If a given impulsive maneuver can significantly change the measurement profile obtained, it improves observability \cite{obsangle, 87}. Hence, the dot product between the ballistic vs forced trajectories (shown in Fig. \ref{obscalc}) is used to measure observability.

Here, $\overline{ \bm{y}}(t)= {\overline{\pmb{r}}^{rel}(t)}/{\parallel \overline{\pmb{r}}^{rel}(t) \parallel}$  denotes the maneuvered measurements of the nominal trajectory and  $\overline{ \bm{y}}_{bal} (t)= {\overline{\pmb{r}}_{bal}^{rel}(t)}/{\parallel \overline{\pmb{r}}_{bal}^{rel}(t) \parallel} $ is the ballistic measurement before a $\Delta \overline{\pmb{v}}$ is applied. 
The dot products are calculated for each segment $i$ to generate the instantaneous observability metric $\eta_i(t)$ as 

\begin{equation}\label{obsindex}
\eta_i(t)=  \overline{ \bm{y}}_{i,bal}(t)^T \overline{ \bm{y}}_i (t) 
\end{equation}
 where $t \in [t_i,t_i + (t_f - t_0)/(n-1)]$ spans the current segment $i$. Note that this follows the intuition that the most observability is gained when the measurement profile fully changes direction, i.e., $\overline{ \bm{y}}_{i,bal}(t)^T \overline{ \bm{y}}_{i}(t) = -1$, and the least observability is present when the measurement is unchanged by the maneuver, 
i.e $\overline{ \bm{y}}_{i,bal}(t)^T \overline{ \bm{y}}_{i}(t) = 1.$  \par 
The objective function contribution of the observability is then obtained as 
 \begin{equation}\label{Gobs}
      G_{\text{obs}} = \sum^{t_f}_{t =0} \eta_i(t) 
  \end{equation}
  where $t_0$ and $t_f$ are the initial and final times of the mission. 

\subsubsection{Safety Metric}\label{SM} 

The safety of the chaser becomes increasingly concerning as it gets closer to the uncooperative target. As observability is optimized, safety constraints must be enforced to prevent the chaser from
coming too close to the target to enhance observability. Two aspects of safety are considered for the design of the nominal trajectory. They are \cite{GG}:
\begin{enumerate}
    \item \textbf{Point-wise (PWS) Safety:} chaser trajectory at time $t_i$ is said to be PWS safe if it is outside a geometrical Keep-Out-Zone (KOZ) defined around the target at time $t_i$ 
    \item \textbf{Passive (PAS) Safety:} chaser trajectory at $t_i$ is said to be PAS safe if it is outside the KOZ at $t_i$ and will remain outside of it for a time $\Delta t$ afterward, during which the chaser will be in uncontrolled flight.  
\end{enumerate}

\begin{figure}[hbt!]
    \centering
    \includegraphics[width=0.8\linewidth]{ 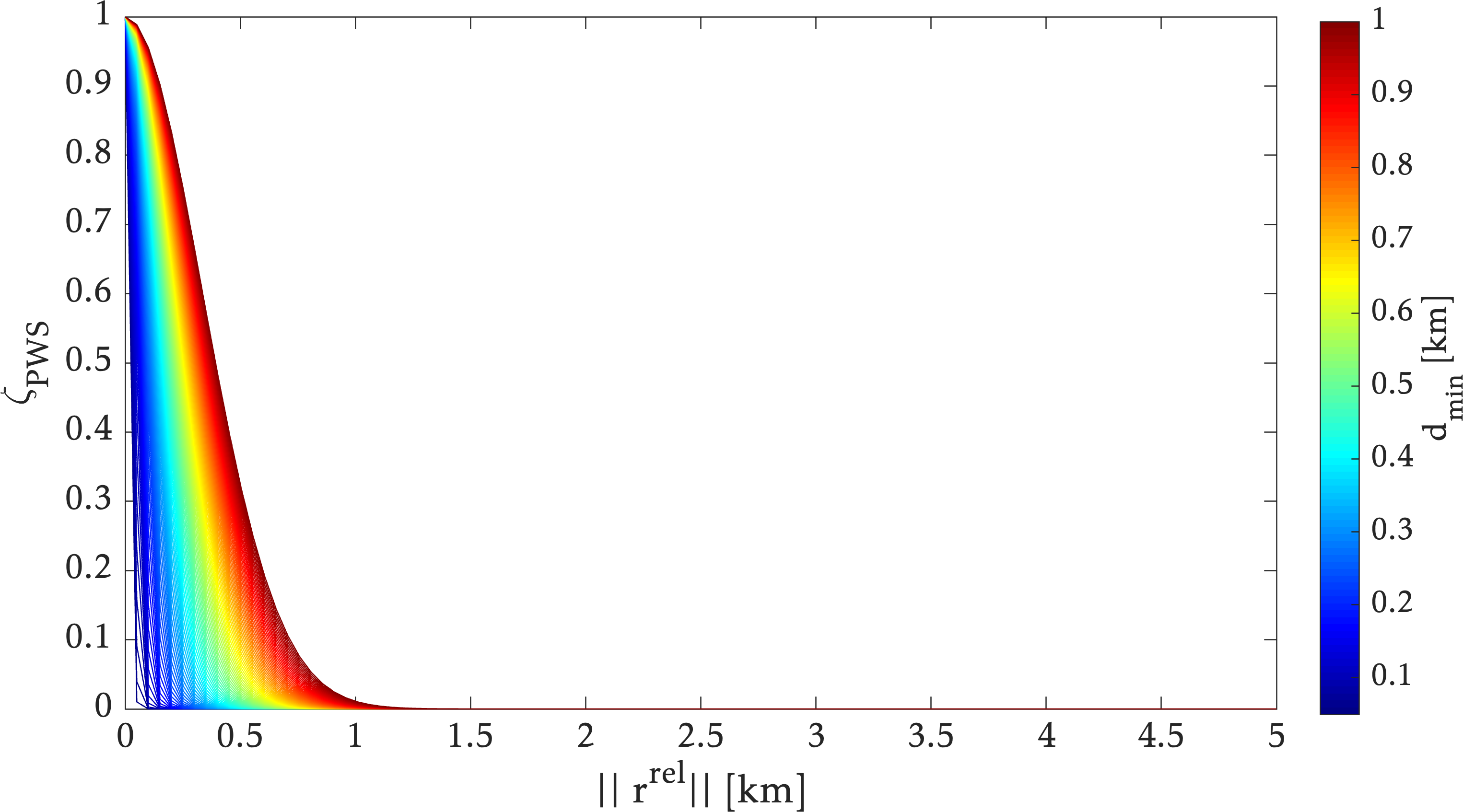}
    \caption{Safety penalty function.}
    \label{safemetricPWS}
\end{figure}

PWS safety can be enforced for this mission by enforcing a distance constraint between the target and the chaser from $t_0$ to $t_f$, such that, 
\begin{equation}\label{PWS1}
    \parallel  \pmb{r}^{rel} (t)\parallel \leq  d_{min}
\end{equation}
where $d_{min}$ is the radius of the KOZ sphere and $ \pmb{r}^{rel}$ is the relative position vector of the chaser with respect to the target. As traditional {PSO} does not allow direct implementation of Eq \eqref{PWS1} as a hard constraint, it is implemented via a Gaussian penalty function as a soft constraint as follows.
\begin{equation}\label{GsafePWS}
 {\zeta_{PWS}}_i(t) =  \exp\left(-\frac{\left(   \parallel  \pmb{r}^{rel}(t) \parallel  \right)^2}{2\sigma^2}\right)
\end{equation}

where $ {\zeta_{PWS}}_{i} (t)$ is the PWS safety index at $t \in [t_i , t_i + (t_f -t_0)/(n-1)]$ and $i$ is the current segment index. In this work,  $\sigma = {d_{min} }/{3}$, allowing $d_{min}$ to be the $3\sigma(99.73\%)$ value. The resultant penalty function is shown in Fig. \ref{safemetricPWS}.

One strategy often explored for PAS safety in rendezvous operations is maintaining a minimum radial-normal (RN) relative distance. This approach is motivated by the higher level of uncertainty present in the tangential relative motion compared to the radial and normal components. It
was first explored by D'Amico in \cite{simonethesis} and has since been utilized in other works \cite{GG, GAIAS2016409}. The minimum RN distance is computed using ROEs. The coordinate conversion from the absolute cartesian states to ROEs can be found in \cite{simonethesis}. The minimum RN distance formulation is given as  

\begin{equation}\label{mindistf1}
    \delta r_{\mathrm{RN}}^{\min }(t)=\min _{t \in[t,t +\Delta t)}   \delta r_{\mathrm{RN}}(t)
\end{equation}

where 
\begin{equation}
      \delta r_{\mathrm{RN}}(t) = \sqrt{(a^T \delta \overline{i} \sin (\overline{u}(t)-\overline{\theta}))^2+(a^T \delta \overline{a}-a^T\delta \overline{e} \cos (\overline{u}(t)-\overline{\theta}-\overline{\phi}))^2}
\end{equation}


where $\overline{a}$ is the semi-major axis of the chaser, $\delta \overline{i} = \overline{i} - i^T$, the difference in inclination between the chaser,  $\delta \overline{a} = \overline{a} - a^T$ is the difference between the semi-major axes,  $\delta \overline{e} = \overline{e} - e^T$  is the difference between the eccentricities,  $\overline{u}(t)$ the mean argument of latitude of the chaser at time $t$, $\overline{\theta}$ the ascending node of the relative orbit and $\overline{\phi}$, the phasing of the relative eccentricity and inclination vectors. A detailed derivation of this distance estimate can be found in \cite{simonethesis}. \par 
Note that Eq. \eqref{mindistf1} is an adequate measure of spacecraft separation only when the spacecraft are in different orbits. When they share 
an orbit, the RN separation distance is always 0. Hence, the true separation distance, including the tangential separation, must be used. \par 

\begin{figure}[hbt!]
    \centering
    \includegraphics[width=0.8\linewidth]{ 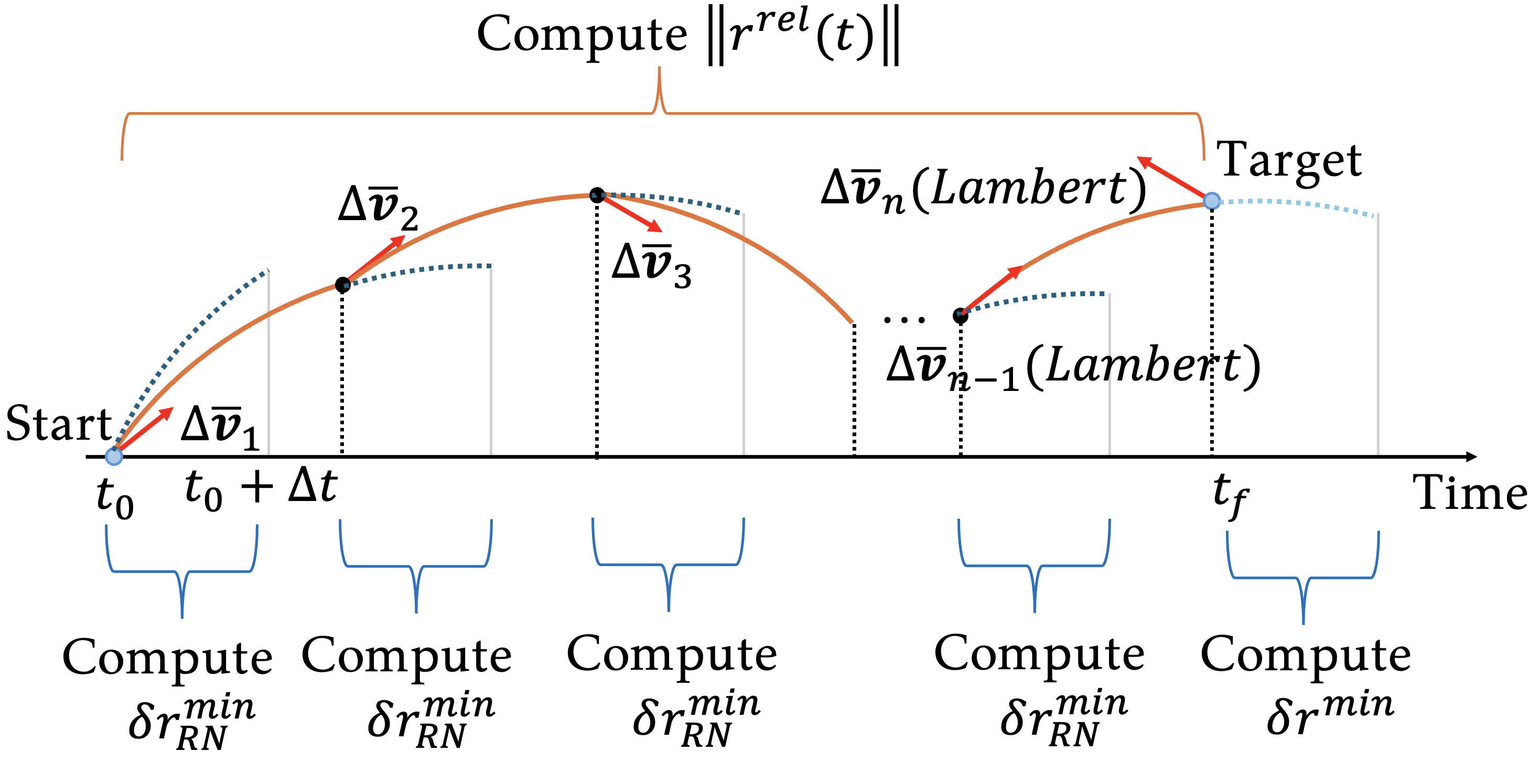}
    \caption{Safety metric: The point-wise safety is calculated on the trajectory (orange), and passive safety is calculated for the resultant ballistic trajectory (dotted blue) if an impulse is missed. }
    \label{safe}
\end{figure}

In this work, the minimum relative distance is calculated for the ballistic segments that would occur in the absence of $\Delta \overline{v}$ applied at the beginning of each segment, as shown in Fig. \ref{safe}. Eq. \eqref{mindistf1} is used to calculate the minimum relative distance at each node up to but not including the last node. At the last step of the trajectory, the chaser and target are in the same orbit; hence, the minimum distance calculation for the last node must consider the total relative distance instead of the separation distance only in the RN plane. 
Therefore, in this work, the minimum distance for PAS Safety is calculated as 

\begin{equation}\label{distPAS}
\delta r ^{\min }_{PAS}(t)=
\begin{cases}
\displaystyle
\text{Eq. } \eqref{mindistf1}, & \quad \text{For nodes } [1, \dots, n-1], \\
\displaystyle
\min _{t \in[t,t +\Delta t)} \parallel \pmb{r}^{rel} (t)\parallel, & \quad \text{For the last node } n.
\end{cases}
\end{equation}

Again, a Gaussian penalty function based on $\delta r ^{\min }_{PAS}$ is used to enforce this safety constraint in the {PSO} as follows. 
\begin{equation}
 {\zeta_{PAS}}_i(t) =  \exp\left(-\frac{\left(\delta r ^{\min }_{PAS}(t)\right)^2}{2\sigma^2}\right)
\end{equation}

where $\sigma = d_{min}/3$. The overall contribution to the {PSO} objective function from safety is then
\begin{equation}\label{Gsafe}
     G_{\text{safety}} =  \sum_{t =t_0}^{t_f} {\zeta_{PAS}}_i(t) +{\zeta_{PWS}}_i(t)
\end{equation}

This calculation of the safety index is summarized in Fig. \ref{safe}. The orange segments denote the actual trajectory, on which Eq \eqref{GsafePWS} is used to determine the PWS safety. The darker dotted blue lines indicate the ballistic trajectory if an impulse is missed, starting from the first node to the second to the last node. For those, the minimum RN plane distance (Eq \eqref{mindistf1}) is used in the PAS safety criterion. The lighter dotted blue line indicates the ballistic arc if the last impulse is missed, for which the total relative distance $\parallel \pmb{r}^{rel} (t)\parallel$ is used in the safety metric, as the target lies in the plane of the target. 

\subsubsection*{Overview}

Algorithm \ref{PSONom} summarizes the steps involved in calculating the objective function $G$ for the nominal trajectory design using {PSO}. Once generated, {PSO} is used to find the optimal nominal $\Delta \overline{\pmb{v}}$ that minimizes $G$. 
\begin{algorithm*}[hbt!!]
\caption{{PSO} objective function calculation.}
\label{PSONom}
 \textbf{Input} $t_0 , t_f, \pmb{x}^T_{0}  ,\overline{\pmb{x}}_0, \overline{\pmb{x}}_f, \Delta \overline{\pmb{v}}$
\begin{algorithmic}
\State Set the chaser state $\overline{\pmb{x}} = \overline{\pmb{x}}_0$.
\State Set the ballistic chaser state $\overline{\pmb{x}}_{bal} = \overline{\pmb{x}}_0$.
\State Set the target state $\pmb{x}^{T} = \pmb{x}^T_0$.
\State Calculate the time between two $\Delta v$ nodes: $\Delta t = (t_f - t_0)/(n-1)$
\State Set $t = 0$.
\For{$i = 1:n-1$} 
\If{$i = n-1$}
\State Calculate $\Delta \overline{\pmb{v}}_{n-1}$ and $\Delta \overline{\pmb{v}}_{n}$ required to reach $\overline{\pmb{x}}_f$ using Lambert's \State method \cite{Vallado2013}.
\EndIf
\State Add the $\Delta v$ to the chaser's state $\overline{\pmb{x}}(4:6) \mathrel{+}=  \Delta \overline{\pmb{v}}_i$.
\State Define a time grid $\delta t= \text{linspace}(t, t+ \Delta t , N_{grid})$. \State  \Comment{$N_{grid}$ is the number of gridpoints in $\delta t$.}
\State Propagate both the maneuvered and ballistic chaser trajectories and \State the target to the next node as shown in Section \ref{propnode}, to get their states \State on $\delta t$.

\State \textit{Safety:} Calculate $\zeta_{\text{PAS},i} (\delta t)$ and $\zeta_{\text{PWS},i} (\delta t)$ for all states on the grid $\delta t$ \State as discussed in Section \ref{SM}
\State \textit{Observability:} Convert the maneuvered and ballistic chaser trajectories \State from absolute ECI to relative as shown in Eq. \eqref{conv1} to obtain $\overline{\pmb{r}}_{bal}^{rel}$ and \State $\overline{\pmb{r}}^{rel}$. Then, calculate $\eta (t)_i$ as discussed in Section \ref{OM}.
\State  Update time $t = t+ \Delta t $.
\State Propagate the chaser and target states using Eq. \eqref{eq21} and set
\begin{equation}
    \overline{\pmb{x}} = \overline{\pmb{x}} (t + \Delta t),\overline{\pmb{x}}_{bal} =\overline{\pmb{x}}_{bal} (t + \Delta t) \ \text{and} \ \pmb{x}^{T} = \pmb{x}^T (t+ \Delta t).
\end{equation}
\EndFor
\State Calculate $G_{\Delta  \overline{\pmb{v}}},  G_{\text{obs}}$ and $ G_{\text{safety}}$ using Eq. \eqref{Gdv}, \eqref{Gobs} and \eqref{Gsafe}, respectively. 
\State Use Eq. \eqref{Gmain} to calculate the combined objective function $G$. 
\end{algorithmic}
\textbf{Output} $G$
\end{algorithm*}

\subsection{Reinforcement Learning-based Guidance Scheme}
As mentioned in the introduction, {RL} is widely used to solve {MDP} problems \cite{montague1999reinforcement}. In fact, many complex high-dimensionality problems can be solved using {RL}, especially in conjunction with Neural Networks \cite{brandonisio2024deep}. In this work, the {RL} algorithm {PPO} \cite{schulman2017proximal} is used due to the continuous nature of the states and actions involved in spacecraft guidance. {PPO} is known to be robust and effective in such environments \cite{schulman2017proximal}, as it can efficiently navigate continuous action spaces by using a policy gradient approach while iteratively adjusting the policy based on sampled trajectories. Furthermore, {PPO} allows clipping of the objective function, which ensures a more stable learning process and aids convergence \cite{HarryThesis}.

\subsubsection{Guidance Problem Formulation}
\begin{figure}[hbt!]
    \centering
    \includegraphics[width=\linewidth]{ 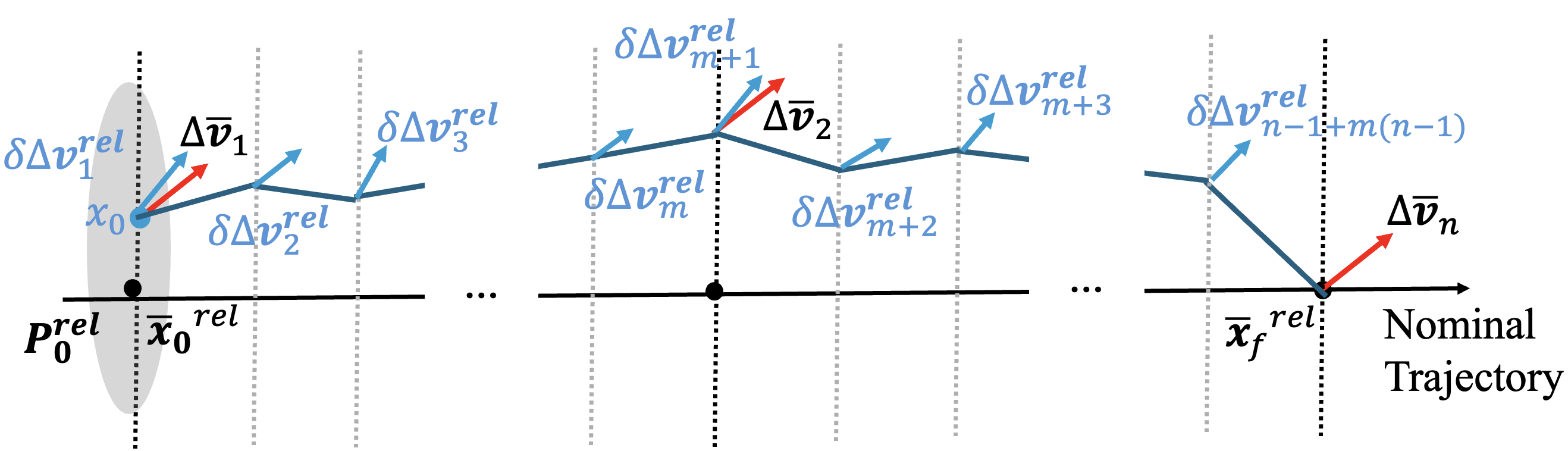}
    \caption{Reinforcement learning guidance.}
    \label{RLprocess}
\end{figure}

The goal of guidance in this work is to retain the observability and safety of the trajectory in the presence of initial state deviations and $\Delta \overline{\pmb{v}}$  errors while still reaching the target state $\pmb{x}_f$/$\pmb{x}_f^{rel}$. A {PPO}-based controller is trained in a stochastic environment to utilize the minimum amount of total $\Delta v$ to accomplish this. Within the stochastic environment, the initial relative state of the chaser is drawn from a uniform distribution of states with mean $\overline{\pmb{x}}^{rel}_0$ and max initial dispersion $\pmb{x}^{rel}_0$, such that

\begin{equation}
    \pmb{x}^{rel}_0 \sim \mathcal{U}(\overline{\pmb{x}}^{rel}_0 - \delta \pmb{x}_i^{\max}, \overline{\pmb{x}}^{rel}_0 + \delta \pmb{x}_i^{\max}).
\end{equation}
Where $\delta \pmb{x}_i^{\max}$ is the maximum initial state dispersion.
Each segment $i$ is allocated $m$ equally spaced additional impulses for guidance. Note that this increases the number of propagation steps from $(n-1)$ to $(n-1) + m(n-1)$ while keeping the overall TOF the same. $j = [1,....(n-1) +m(n-1)]$ is used to iterate through these sub-segments. The additional $\Delta v$ utilized for guidance is again denoted $\delta \Delta \pmb{v}^{rel}$. This formulation is illustrated in Fig. \ref{RLprocess}.
  
As mentioned, the primary goal of the guidance scheme is to reach the target state $\pmb{x}_f$ at $t_f$ while consuming minimal additional fuel without compromising safety or observability. Hence, 
in this guidance scheme, the spacecraft is allowed to deviate from its {PSO}-generated nominal trajectory to some degree to optimize its objectives. At time $t_j$, such a deviation along the trajectory is denoted by  
\begin{equation}
    \delta \pmb{x}_j = \pmb{x}^{rel}_j - \overline{\pmb{x}}^{rel}_j.
\end{equation}
This deviation can be linearly propagated in time to obtain the next state $\delta \pmb{x}_{j+1}$ using the Yamanaka-Ankersen STM \cite{doi:10.2514/2.4875} such that 
\begin{equation}\label{qcqpc4}
\delta \pmb{x}_{j+1} = \bm{\Phi}_{{YA}_j} \delta \pmb{x}_j  +  \bm{\Phi}_{{YA}_j} M \delta \Delta \pmb{v}_j^{rel}.
\end{equation}
The propagated state ($\delta \pmb{x}_{j+1}$) can then be constrained based on the previous state as 
\begin{equation}\label{qcqp3}
 \delta \pmb{x}_{j+1}^T  \delta \pmb{x}_{j+1} \leq {\alpha}^2   \delta \pmb{x}_j^T   \delta \pmb{x}_j. 
\end{equation}
Here, $\alpha$ is a contraction constant that determines the degree of deviation. In other words, adjusting $\alpha$ influences how close or far the guided trajectory will be from the nominal trajectory at the next time step. The $\delta \Delta \pmb{v}_j^{rel}$ - the $\Delta v$ required to reach the new state that abides by the constraint given in Eq. \eqref{qcqp3}- can then be found by solving the following {QCQP} problem. 

\begin{align}\label{qcqpRL}
   \min_{\delta \Delta \pmb{v}^{rel}_j}  \quad &   \parallel\delta \Delta \pmb{v}^{rel}_j  \parallel \\
    \text{subject to} \quad & \text{Eq. \eqref{qcqp3}}\\
                             & \text{Eq. \eqref{qcqpc4}}
\end{align}

This can be solved via a quadradic solver. Varying $\alpha$ results in variations in the inequality constraint presented in Eq. \eqref{qcqp3}, resulting in trajectories with differing $\Delta v$ consumptions, observability, and safety characteristics. This work uses {PPO} to train a neural network that selects $\alpha$ values for the whole trajectory to minimize the fuel cost and maintain safety and observability.




\subsubsection{PPO Objective Function} 
At each step $j$, the objective function $R_j$  (known as the reward function) that the {PPO} must maximize when selecting $\alpha$ is given by 
\begin{equation}\label{reward}
    R_j = 
    \begin{cases} 
        -\left( \delta \Delta \pmb{v}^{rel}_j + P_{j,safe} + P_{j,obs} \right) & \text{if } j \bmod (m+1) \neq 0, \\
        -\left( \delta \Delta \pmb{v}^{rel}_j + \Delta \overline{\pmb{v}}^{rel}_{\left\lfloor \frac{j}{m+1} \right\rfloor} + P_{j,safe} + P_{j,obs} \right) & \text{otherwise.}
    \end{cases}
\end{equation}

Note that this allows the total fuel consumption to be minimized, as it minimizes $\ \delta \Delta \pmb{v}^{rel}_j + \Delta \overline{\pmb{v}}^{rel}_{\left\lfloor \frac{j}{m+1} \right\rfloor}  \  $ when the guidance node $j$ coincides with a node where a nominal $\Delta \overline{v}$ is applied. Note that ${\left\lfloor \Box \right\rfloor}$ is the floor operator. $P_{j, safe}$ denotes a penalty applied to prevent unsafe trajectories. This is done by calculating $\parallel \pmb{r}^{rel}_j\parallel$ and $\delta r ^{\min }_{j,PAS}$ as shown in Eq. \eqref{distPAS}. Then, the safety penalty is imposed as 
\begin{equation}\label{Psafe}
P_{j,safe} =  \begin{cases} 
 \rho_{safe}(d_{min}- \parallel \pmb{r}^{rel}_j\parallel) & \text{if }   \parallel \pmb{r}^{rel}_j\parallel \leq d_{min} \\
  \rho_{safe} (d_{min} - \delta r ^{\min }_{j,PAS} )& \text{else if }   \delta r ^{\min }_{j,PAS} \leq d_{min} \\
 0 & \text{otherwise } 
\end{cases}
 \end{equation}
 where $\rho_{safe}$ is a user-imposed constant penalty value. $P_{j, obs}$ denotes a penalty applied to prevent unobservable trajectories. This is done by calculating the instantaneous observability metric $\eta_j$ via Eq. \eqref{obsindex} and comparing it to the $\eta$ of the nominal trajectory at $t_j$ as follows. 

\begin{equation}\label{Pobs}
P_{j,obs} =  \begin{cases} 
 \rho_{obs} (\eta_j - \bar{\eta}(t_j)) & \text{if }   \eta_j \geq \bar{\eta}(t_j)\\
 0 & \text{otherwise } 
\end{cases}
 \end{equation}
 where $\rho_{obs}$ is a user-imposed constant penalty value and $\overline{\eta}$ is the observability metric value from the nominal trajectory.
Note that this only imposes a penalty on observability if it is reported to be worse than the nominal trajectory following a maneuver. While it does not check for observability between sub-segments, it was seen to be sufficiently adequate as long as the sub-segments are not placed too far apart. 

\subsubsection{PPO Actor Input}
In {PPO}, the RL actor must be provided with a set of inputs that contain enough information to estimate the potential reward of an action. This work uses $\bm{S}$ to denote the input state space provided to the agent at each propagation step $j$. The parameters to include in $\bm{S}$ were determined by trial and error. The best training outcomes were obtained by including the remaining number of nominal $\Delta \overline{v}$ nodes $n_{rem}$, the ECI state of the chaser $\pmb{x}_j$, the relative state of the chaser $\pmb{x}_j^{rel}$, previous reward $R_{j-1}$ and the previous $\alpha_{j-1}$ and the difference between the current relative state and the relative nominal state at time step $j$ in $\bm{S}$.
Hence: 
\begin{equation}
    \bm{S}_j = [n_{rem}, \pmb{x}_j , \pmb{x}^{rel}_j, R_{j-1}, \alpha_{j-1},\pmb{x}^{rel}_j- \overline{\pmb{x}}^{rel}_j ]
.\end{equation}
As it is advisable to normalize inputs to ensure stability and consistent feature scaling \cite{fereoli2022meta}, a simple min-max  normalization \cite{lafarge2023adaptive} is used to normalize $\bm{S}$
such that $\bm{S} \in [-1,+1]$ such that 
\begin{equation}\label{actorin}
\bm{S}_j^*=2\left(\frac{\bm{S}_j-\bm{S}_{\min}}{\bm{S}_{m a x}-\bm{S}_{\min}}\right)-1.
\end{equation}

\subsubsection{PPO Actor Output}  
As mentioned, the {PPO} outputs the parameter $\alpha$ at each iteration. 

 \subsubsection*{Overview} 
The process inside the {PPO} guidance environment for a single step $j$ is given in Algorithm \ref{GE}. Note that $\overline{\pmb{x}}^{rel}_{j =[ 1 ,..., n-1 + m(n-1)}$ indicates the nominal spacecraft trajectory from the {PSO} optimization. At the start of the guidance environment, 
\begin{equation}
    \pmb{x}_{j = 0} =  \pmb{x}(t_0) ,   \pmb{x}^T_{j = 0} =  \pmb{x}^T(t_0) \ \text{and} \ n_{rem} =  n-1 + m(n-1)
\end{equation}
where $\pmb{x}(t_0) $ is derived by converting $ \pmb{x}^{rel}_0$ to absolute ECI coordinates.

\begin{algorithm*}[hbt!!]
\caption{Guidance Environment.}
\label{GE}
 \textbf{Input} $\alpha_j , t_0 , t_f ,\pmb{x}_j ,\pmb{x}_{j,bal} , \pmb{x}^T_j, \Delta \overline{\pmb{v}}, n , m , \overline{\pmb{x}}^{rel}_j , \bar{\pmb{x}}^{rel}_{j+1}$ and $\overline{\eta} _j$
\begin{algorithmic}
\State Calculate time at $j+1$ as $t_{j+1} = t_j + {(t_f - t_0)}/{(n-1 + m(n-1))}$.
\If{$j \bmod (m+1) == 0 $} \State \Comment{If $j$ is a node with a nominal $\Delta \overline{\pmb{v}}$ implementation, add it to the state.}
\State $\pmb{x} = \pmb{x} + [\bm{0}_{3 \times 1} ;\Delta \overline{\pmb{v}}_{{\left \lfloor {j/(m+1)}  \right \rfloor)} } ]$
\State Convert $\Delta \overline{\pmb{v}}_{\left \lfloor {j/(m+1)}  \right \rfloor}$ to RTN to obtain  $\Delta \overline{\pmb{v}}^{rel}_{\left \lfloor {j/(m+1)}  \right \rfloor}$
\EndIf 
\State Convert $\pmb{x}_j$ to the relative state $\pmb{x}^{rel}_{j}$.
\State Calculate the deviation of the current state from the nominal. $\delta \pmb{x}^{rel}{_j} =\pmb{x}^{rel}_{j} - \overline{\pmb{x}}^{rel}_{j}$.
\State \textbf{if} $\eta_{rem} = 1$ \textbf{then} $\alpha = 0$ \Comment{This is to ensure that the target $\pmb{x}_f$ is reached.}
 \State Find $\delta \Delta \pmb{v}^{rel}_j$ by solving the QCQP in Eq. \eqref{qcqpRL}.
\State Convert $\delta \Delta \pmb{v}^{rel}_j$ to an absolute ECI quantity, defined as  $\delta \Delta \pmb{v}_j$.
\State Update the state as:  $\pmb{x} = \pmb{x} + [\bm{0}_{3 \times 1} ;\delta \Delta \pmb{v}_j]$.
\State Propagate both the maneuvered and ballistic chaser trajectories and the target to $t_{j+1}$.
\State \textit{Observability:} Calculate $P_{j+1,obs}$ as shown in Eq. \eqref{Pobs}.
\State \textit{Safety:} Calculate $P_{j+1,safe}$ as shown in Eq. \eqref{Psafe}.
\State Update $n_{rem} = n_{rem} -1$.
\If{$n_{rem} == 0$} 
\State {\textbf{Stop.}}
\Else 
\State {Update states and time: $t_j = t_{j+1}, \pmb{x}^T_j =  \pmb{x}^T_{j+1},\pmb{x}_j =  \pmb{x}_{j+1} $ and $\pmb{x}_{j,bal} =  \pmb{x}_{j+1,bal}$}
\State {\textbf{Repeat from the start.}}
\EndIf
\State Calculate the Reward $R_{j+1}$ using Eq. \eqref{reward} and calculate the actor input state  $\bm{S}_{j+1}^*$ using Eq.  \eqref{actorin}. 
\end{algorithmic}
\textbf{Output} $\bm{S}_{j+1}^*, R_{j+1}$ 
\end{algorithm*}

\section{Numerical Simulation Results}

\subsection{Simulation Conditions}
The simulation parameters are provided in Table \ref{coordinates} and were set to recreate the test case provided in \cite{YUAN2022107812}, which performs a V-bar far-range transfer. The relevant initial state was set as $\pmb{x}_{0}^{rel} = [\SI{1}{\meter}, \SI{100}{\kilo\meter}, \SI{5}{\meter}, \SI{-0.02}{\meter\per\second},  \SI{0.01}{\meter\per\second},0]^T$ and the final state was defined to be  $\pmb{x}_{f}^{rel} =[0, \SI{1}{\kilo\meter},  \SI{5}{\meter}, 0,  0,0]^T $.  The rendezvous target was set to be an object in a circular orbit of 300 km altitude, with 99.8 deg inclination.  The Keplerian coordinates of the target at the initial time was $\pmb{\text{\oe}}^T (t_0) = [\SI{300}{\kilo\meter}+R_\oplus, 0, \SI{99.8}{\degree},0,0,0]^T$, where $R_\oplus = \SI{6378,136}{\kilo\meter}$ is the radius of the Earth. The time allocated for the far-range approach (TOF) was equivalent to four orbital periods of the target, totaling 6 hours.  For the safety criteria, the KOZ radius $d_{min} $ was set to 500 m, and the duration $\Delta t$ for the PAS safety evaluation was set to be one orbital period of the target, similar to \cite{GG}. 

\subsection{Extended Kalman Filter Parameters}\label{EKFdisc}
 In this work, an EKF was utilized to verify the observability of the nominal trajectory and the RL-guided trajectories by calculating the maximum eigenvalue of the position covariance. During this process, the initial estimate of the state $\hat{\pmb{x}}^{rel}$ for the EKF was set to be 
\begin{equation}
    \hat{\pmb{x}}^{rel} = \overline{\pmb{x}}^{rel} + \Delta \overline{\pmb{x}}^{rel} \ \text{where} \ \Delta \overline{\pmb{x}}^{rel} \sim \mathcal{N}(\bm{0}_{6 \times 1}, \sigma_{\Delta \overline{\pmb{x}}^{rel}})
\end{equation}
 and $\sigma_{\Delta \overline{\pmb{x}}^{rel}} = [\SI{300}{\meter}, \SI{100}{\meter}, \SI{100}{\meter} ,\SI{0.3}{\meter\per\second}, \SI{0.3}{\meter\per\second} , \SI{0.3}{\meter\per\second}]^T$ set to match \cite{doi:10.2514/1.G00683}. The covariance of the EKF was reset to $\text{diag}(\sigma^2_{\Delta \overline{\pmb{x}}^{rel}})$ at each impulse. The EKF time step was set to 10 s.  The sensor noise covariance was established as $\bm{\sigma}_s = [\SI{1e-3} ,\SI{1e-3}, \SI{1e-3}]^T$ while the process noise covariance was set as $\bm{\sigma}_w =  [\bm{\sigma}_{w_r}, \bm{\sigma}_{w_v}]^T = [\SI{50}{\meter},\SI{50}{\meter},\SI{50}{\meter},\SI{0.1}{\meter\per\second},\SI{0.1}{\meter\per\second},\SI{0.1}{\meter\per\second}]^T$
corresponding to the state and sensor errors provided in \cite{YUAN2022107812}. Consequently, the observations for the EKF were obtained as 
 \begin{equation}
     \pmb{y}(t) = \frac{ {\pmb{r}}^{rel}_E(t) }{\parallel {\pmb{r}}^{rel}_E(t)\parallel } + \delta   \pmb{y}_i(t)  \ \text{where}  \  \pmb{y}_i(t)  \sim \mathcal{N}(\bm{0}_{3 \times 1}, \bm{\sigma}_s)
 \end{equation}
 where 
 \begin{equation}
     \pmb{r}^{rel}(t)_E=   \pmb{r}^{rel}(t) + \delta  \pmb{r}^{rel}(t) \ \text{where} \ \delta  \pmb{r}^{rel}  \sim \mathcal{N}(\bm{0}_{3 \times 1}, \bm{\sigma}_{w_r}).
 \end{equation}
 
 These parameters are also provided in Table \ref{coordinates}. The formulation of the EKF filter used in this work can be found in \cite{franquiz2019trajectory}.

\begin{table*}[hbt!]
\centering 
\begin{tabular}{lll} \hline 
Parameter                   & Value                          & Units     \\\hline 
\multicolumn{3}{l}{\textbf{Target spacecraft environment}}                                 \\
Target altitude                    & 300                            & km        \\
Target orbital period, T           & 90                             & min       \\
Initial Keplerian target state  $\pmb{\text{\oe}}^T (t_0)$             & $[\SI{300}{\kilo\meter}+R_\oplus, 0, \SI{99.8}{\degree},0,0,0]^T$  & km deg     \\  
\textbf{Far-range mission parameters} & \textbf{}                      & \textbf{} \\
Initial state  $\pmb{x}_{0}^{rel}$             & $[1, 100,000, 5,-0.02,0.01,0]^T$ & m m/s     \\
Final state     $\pmb{x}_{f}^{rel}$            & $[0,1000,0,0,0,0]^T$         & m m/s     \\
Rendezvous duration, $TOF$     & 4 (6)                             & Orbits (hours)   \\ 
KOZ radius                  & 500                            & m         \\
PAS duration                & 1                              & Orbit     \\  
\textbf{EKF parameters}         & \textbf{}                      & \textbf{} \\
Initial errors, $\sigma_{\Delta \overline{\pmb{x}}^{rel}}$          & {[}300 100 100 0.3 0.3 0.3{]}     & m m/s     \\
Navigation time step        & 10                             & s     \\ 
Sensor noise covariance   $\sigma_s$     & $\SI{1e-3}{}$                            & rad /axis     \\ 
Process noise covariance $1 \sigma_w$      & $[50, 50, 50, 0.1, 0.1, 0.1]^T$                           & m m/s     \\  
\textbf{PSO parameters}         & \textbf{} \\
Number of $\Delta \overline{v}$ impulses          & 6                            &           \\
Maximum impulse magnitude $\Delta \overline{v}_{lim}$  & 10                             & m/s       \\ 
$N_{grid}$ in Algorithm \ref{PSONom} & 1000 & \\ 
\textbf{RL parameters}         & \textbf{} \\
Number of $\delta \Delta v$ impulses & 34  & \\ 
RL initial condition dispersion $\delta \pmb{x}_i^{\max}$ & $[15,15,0.75,0,0,0]^T$ &km m/s\\ 
Observation penalty $\rho_{obs}$ & 1000 & \\ 
Safety penalty $\rho_{safe}$ & 1000 & \\ 
\hline    
\end{tabular}
\caption{Problem Parameters (extracted from \cite{YUAN2022107812} and \cite{doi:10.2514/1.G000822}).}
\label{coordinates}
\end{table*}


6 nominal $\Delta \overline{v}$ nodes were utilized when obtaining the nominal trajectory for this test case, equally spaced in time. Their maximum impulse $\Delta \overline{{v}}_{\text{lim}}$ was set to be $\SI{10}{\meter\per\second}$ to match the condition set in \cite{YUAN2022107812}. \par
\begin{figure*}[hbt!!]
    \centering
    \includegraphics[width=1\linewidth]{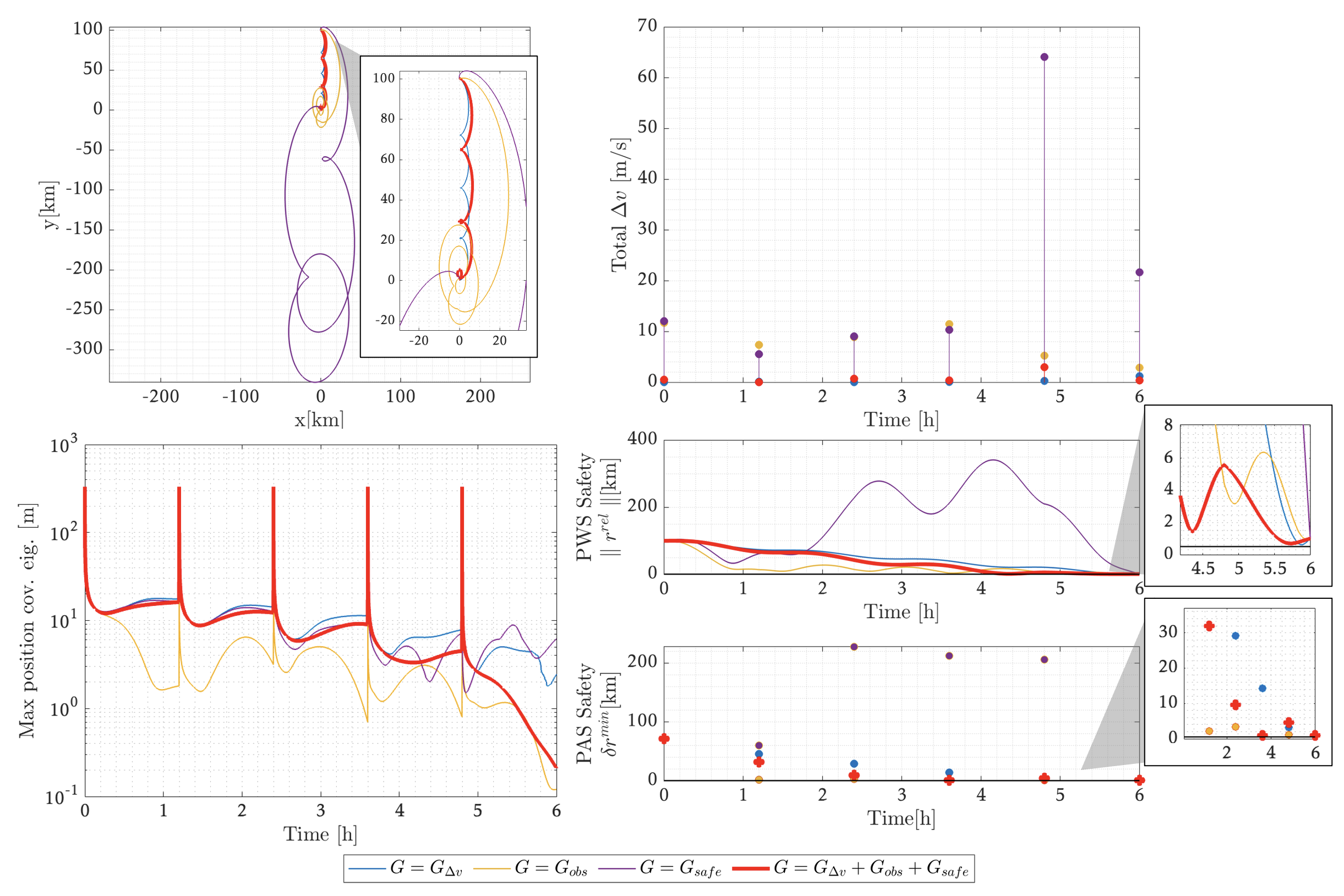}
    \caption{4 h nominal trajectory. Top left: relative trajectories. Top right: Total $\Delta v$ consumption. Bottom left: Max position covariance observed via Extended Kalman Filter (observability). Bottom Right: Relative distance (Safety).}
    \label{refbig}
\end{figure*}

The nominal trajectory obtained is shown in Figure \ref{refbig}, and the corresponding absolute ECI $\Delta \overline{v}$s imparted by the chaser is given in Table \ref{psodv}. It is shown alongside the trajectories derived by optimizing the individual cost functions. The generated nominal trajectory consumes
 \SI{5.0381}{\meter\per\second} $\Delta v$ in total, which is \SI{0.44}{\meter\per\second} less than the $\Delta v$ consumption given in \cite{YUAN2022107812} under these conditions. This difference is likely due to the approach cone constraint utilized to maintain PWS safety in \cite{YUAN2022107812}. In this work, PAS and PWS safety are considered and maintained by establishing a KOZ. This safety strategy appears more optimal than an approach cone regarding fuel consumption.
 \par 
As mentioned, the estimated state uncertainty from AON is measured by the maximum eigenvalue of the position covariance obtained from the EKF discussed in Section \ref{EKFdisc}. The maximum position covariance is \SI{0.209}{\meter} at the end of the nominal trajectory, which illustrates good observability. As a percentage of range, this navigation error is 0.0209\%, which is smaller compared to the 0.73\% error reported in \cite{YUAN2022107812}.
 \par
 The minimum PAS safety distance calculated as shown in Equation \eqref{distPAS} and the distance between the Servicer and the debris throughout the mission is also shown in Figure \ref{refbig}. It can be seen that the nominal trajectory and the PAS relative distances remain outside the KOZ zone, even though the KOZ boundary is reached at some points. Hence, the nominal trajectory is PAS and PWS safe. \par 

\begin{table}[hbt!!]
\centering 
\begin{tabular}{cc}
\hline
Node No. $n$                                                & $\Delta \overline{\pmb{v}}$ (m/s) \\ \hline
1                                                           & $ [-0.3222  , -0.0691  ,  0.3964]$   \\
2                                                           &  $[ -0.0000  ,  0.0000 ,  -0.0000]$   \\
3                                                           &   $[  -0.7197   , 0.0313 ,  -0.1560]$    \\
4                                                           &    $ [ 0.2889 ,  -0.0495 ,   0.2627]$    \\
5                                                          &    $ [ -0.0488 ,   0.5096 ,  -2.9538]$    \\
6                                                           &    $ [ -0.2966 ,   0.0456  , -0.2554]$    \\ \hline 
$\Sigma\parallel \Delta \overline{\pmb{v}}\parallel $ (m/s) & 5.0381                          \\ \hline
\end{tabular}
\caption{Optimized nominal $\Delta \bar{v}$ values.}
\label{psodv}
\end{table}

\subsection{RL Guidance}
For the RL guidance, initial state of dispersion was set to be $\delta \pmb{x}_i^{\max} = [\SI{15}{\kilo\meter},\SI{15}{\kilo\meter},\SI{0.75}{\kilo\meter}, 0, 0,0]^T $ to match the test case given in \cite{YUAN2022107812}. The number of $\delta \Delta v$ impulses introduced was set to 34, so the total number of nodes can be 40, same as in \cite{YUAN2022107812}. The penalty values for observability and safety $\rho_{obs}$ and $\rho_{safety}$ were set to 1000. $\alpha_{max} $ was set to $2$ by trial and error, as limiting the contraction values below or above led to suboptimal outcomes.  \par 
This guidance scheme was developed using the PPO algorithm from Stable-Baselines 3 \cite{raffin2021stable}. The hyperparameters of the PPO chosen for this work are listed in Table \ref{hypers}, which were set by trial and error. These play a crucial role in the behavior and performance of the RL agent during its training and can affect the training process and the final results. The actor and critic neural networks were assigned four hidden layers and an initial standard deviation of 0.25.
\begin{table}[hbt!]
\centering
\begin{tabular}{ l c }
\hline
{Parameter}                         & {Value} \\  \hline 
Batch size                        & 64    \\
Number of epochs                  & 10    \\
Number of episodes for evaluation & 6    \\
GAE lambda $\lambda$               & 1     \\
Discount factor $\gamma$          & 0.99  \\
Clip parameter $\epsilon$         & 0.1   \\
Learning rate $\alpha$            & 0.003 \\
Entropy coefficient               & 0.01  \\ \hline
\end{tabular}
\caption{Selected hyperparameter values. The definitions of these hyperparameters can be found in \cite{schulman2017proximal}.}
\label{hypers}
\end{table}

\subsubsection{Training}
When training the RL guidance, the contraction parameter $\alpha$ was scaled to $\alpha \in [-1,1]$ outside the environment and rescaled to be inside the range of $[0,\alpha_{max}]$ in the environment. This scaling was done to ensure symmetry in the neural network, which was seen to produce better results.

While the RL training was conducted with stochastic initial states, the policies were evaluated on a sparse grid of points from the initial uniform distribution. These points were computed as 
\begin{equation}\label{sigpointCalc}
\begin{array}{ll}
\bm{s}_k=\overline{\pmb{x}}^{rel}_{0} +0.5\delta \pmb{x}_{i,k}^{\max} , & \forall i \in\{1, \ldots, 6\}, \\
\bm{s}_k=\overline{\pmb{x}}^{rel}_{0}-0.5\delta \pmb{x}_{i,k}^{\max}, & \forall k \in\{7, \ldots, 12\} .
\end{array}
\end{equation}This strategy for policy evaluation was employed to ensure equity when testing the developed policies. The best policy was then selected based on the performance on the $\bm{s}$ points. 




\subsubsection{Testing}
Following training, the best PPO policy was tested using an MC campaign with 500 different initial states within the covariance $\bm{P}_0^{rel}$, along with errors in the nominal $\Delta \overline{\pmb{v}}$ s. Its performance was then compared against several other benchmarking $\alpha $ strategies. This was done to measure the performance of the policy statistically and to evaluate the robustness of the guidance scheme developed.

\textbf{Thrust errors} \par 

Note that the errors considered here are slightly higher than those considered in \cite{YUAN2022107812}, but they are only applied to the nominal $\Delta v$ nodes.

\textbf{Benchmark $\alpha$ Strategies} \par 

Several $\alpha$ strategies were also developed to evaluate the performance of the RL-generated $\alpha$ ($\alpha_{RL}$) in comparison. These distributions include:
\begin{itemize}
    \item {Linearly decreasing $\alpha$  ($\alpha_{LD}$):} where $\alpha$ linearly decreases from 1 to 0 as time goes from $t_0$ to $t_f$, hence
\begin{equation}
    \alpha_{{LD}_j} = 1 - \frac{t_j - t_0}{t_f - t_0}.
\end{equation}
Note that this allows consistently decreasing deviations from the nominal trajectory and ensures that at the last node $\alpha_{LD} = 0$.
\item {Constant $\alpha$ ($\alpha_C$):} where $\alpha$ is consistently 0. 
\begin{equation}
    \alpha_{C}= 0
\end{equation}
Note that this requires the trajectory to immediately reach the nominal path upon reaching $j =2$. This ensures that the trajectory would have the safety and observability characteristics of the nominal path almost from the start. 

\item  Midpoint-optimized $\alpha$ ($\alpha_{S}$): where $\alpha$ is obtained by optimizing the total reward obtained by the $\bm{s}$ points obtained from Eq \eqref{sigpointCalc}. 
I.e,  $\alpha_{S}$ is obtained by solving 
\begin{equation} \label{sigpoint}
   \alpha_S = \arg \min_{\alpha} \left( \sum_{k=1}^{12} \sum_{j = 1}^{j_{\text{end}}} R_j\left(\pmb{x}_0^{rel} = \bm{s}_{k}, \alpha = \alpha_j\right) \right).
\end{equation}
where $R_j$ is the reward shown in Eq. \eqref{reward} and $j_{end} = n-1 + m(n-1)$. This optimization was done via interior-point optimization \cite{Byrd2000ATR}.
\end{itemize}

\subsubsection{Results}
\begin{figure}[hbt!]
    \centering
    \includegraphics[width=0.8\linewidth]{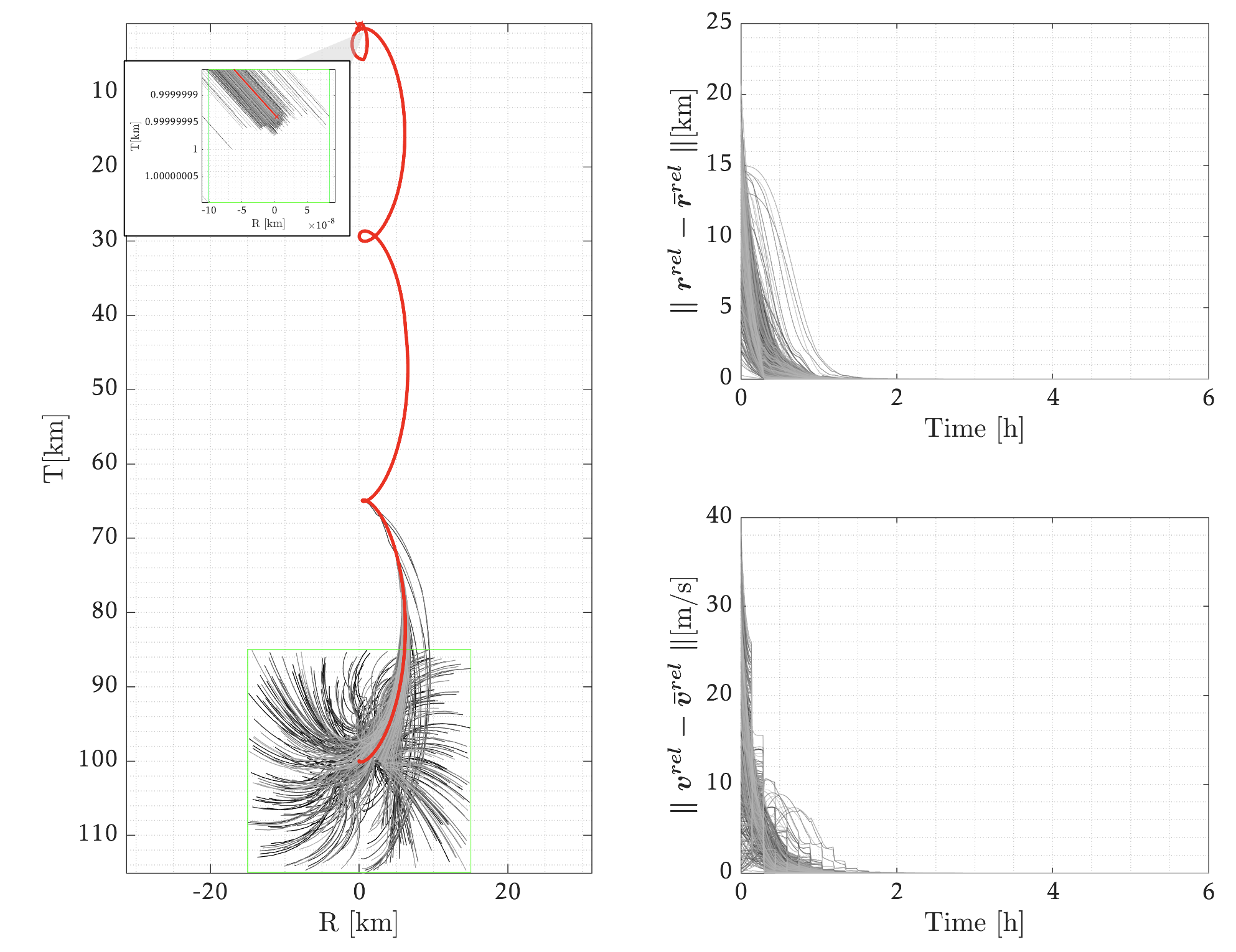}
    \caption{Chaser trajectories from the MC simulation (no error). Dotted lines denote where impulses are applied. The thick blue line is the nominal trajectory.}
    \label{trajlong}
\end{figure}

\begin{figure}[hbt!]
    \centering
    \includegraphics[width=0.8\linewidth]{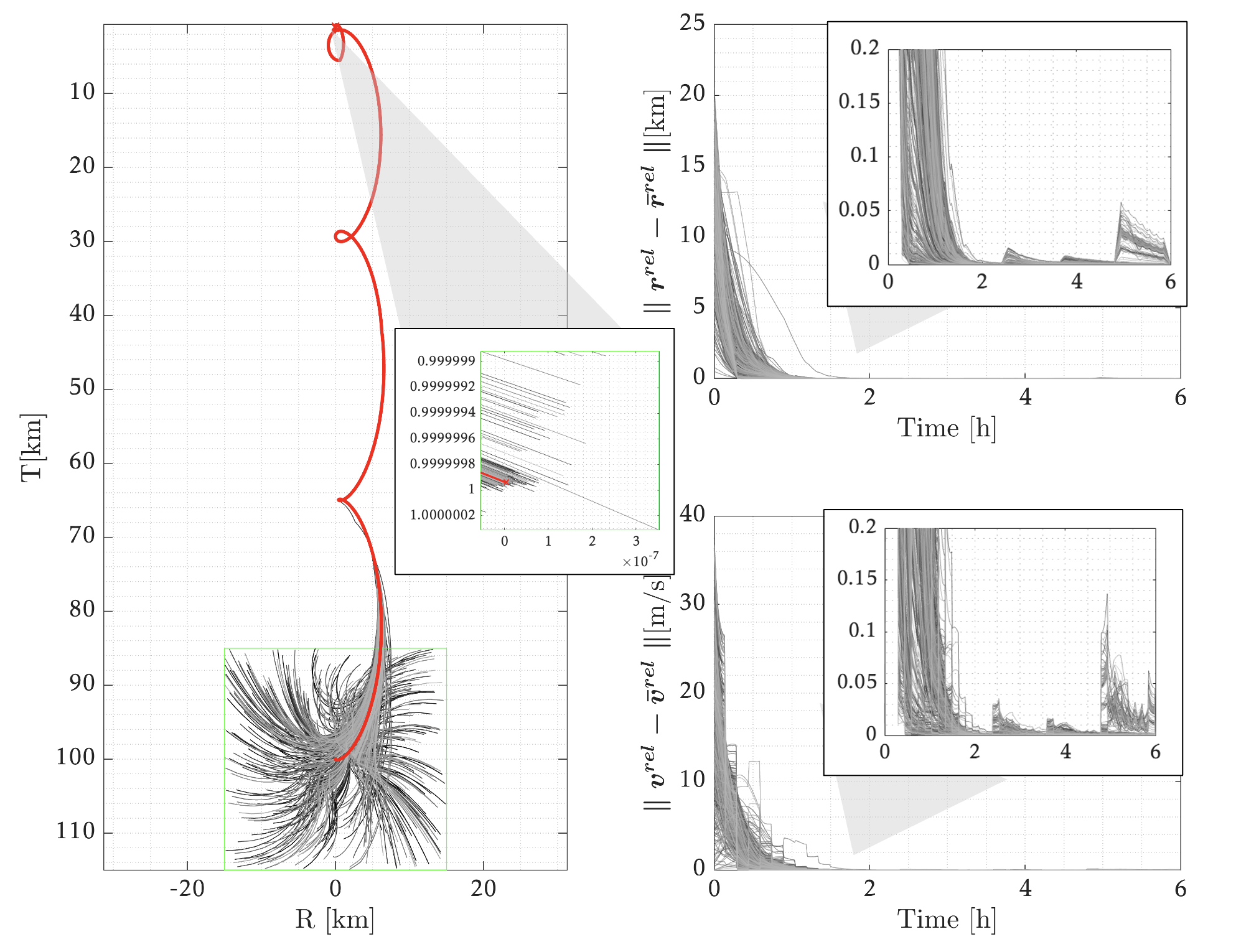}
    \caption{Chaser trajectories from the MC simulation (low error). Dotted lines denote where impulses are applied. The thick blue line is the nominal trajectory.}
    \label{trajlong1}
\end{figure}

\begin{figure}[hbt!]
    \centering
    \includegraphics[width=0.8\linewidth]{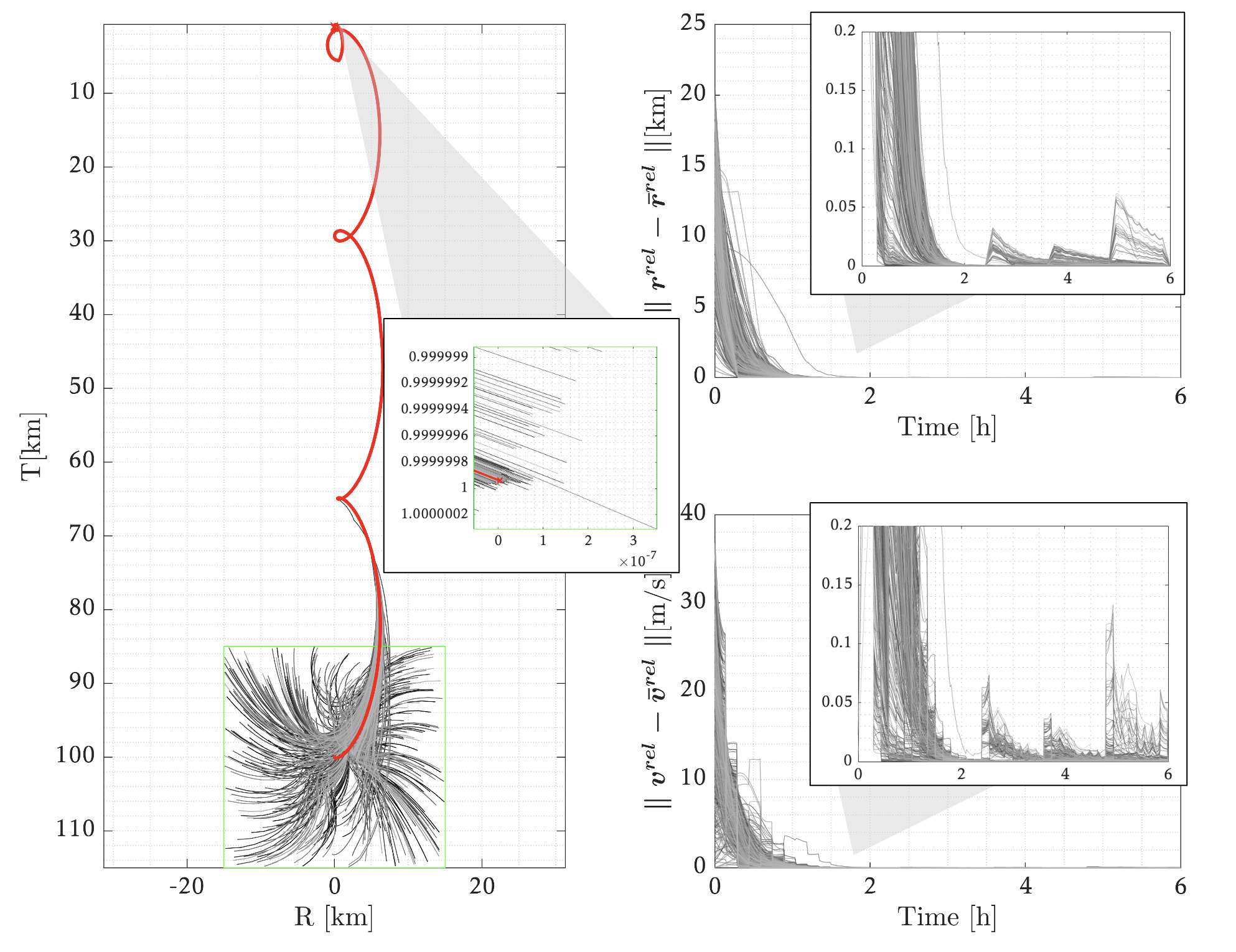}
    \caption{Chaser trajectories from the MC simulation (higherror). Dotted lines denote where impulses are applied. The thick blue line is the nominal trajectory.}
    \label{trajlong2}
\end{figure}

 Fig. \ref{trajlong}, \ref{trajlong1} and \ref{trajlong2} show the trajectories obtained via the 500 sample- MC simulations for the case without any $\Delta \overline{v}$ errors, with low errors and with high errors, respectively. Notably from the no error case it can be seen that the RL-trained $\alpha$ reaches the nominal trajectory at approximately 2 h from $t_0$, a significant time before the end of the transfer. The low error and high error case shows spikes in position and velocity deviation from the nominal trajecotry due to the $\Delta \overline{v}$ errors introduced.

\begin{figure}[hbt!]
    \centering
    \includegraphics[width=\linewidth]{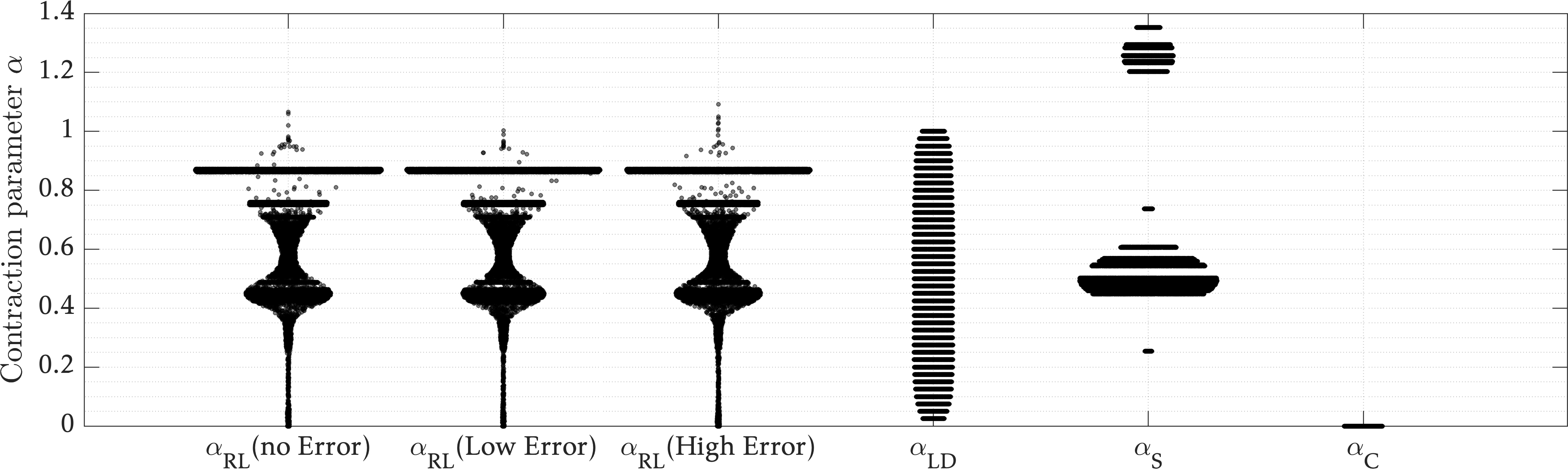}
    \caption{Distribution of $\alpha$ for RL-optimized vs $\alpha$ obtained by other means. }
    \label{alphalong}
\end{figure}

The $\Delta v$ consumption, observability, and robustness of $\alpha_{RL}$ is compared against that of $\alpha_{LD}$, $\alpha_C$ and $\alpha_S$. The corresponding distributions of $\alpha$ are shown in Fig. \ref{alphalong}. \par

\begin{figure}[hbt!]
    \centering
     \includegraphics[width=0.8\linewidth]{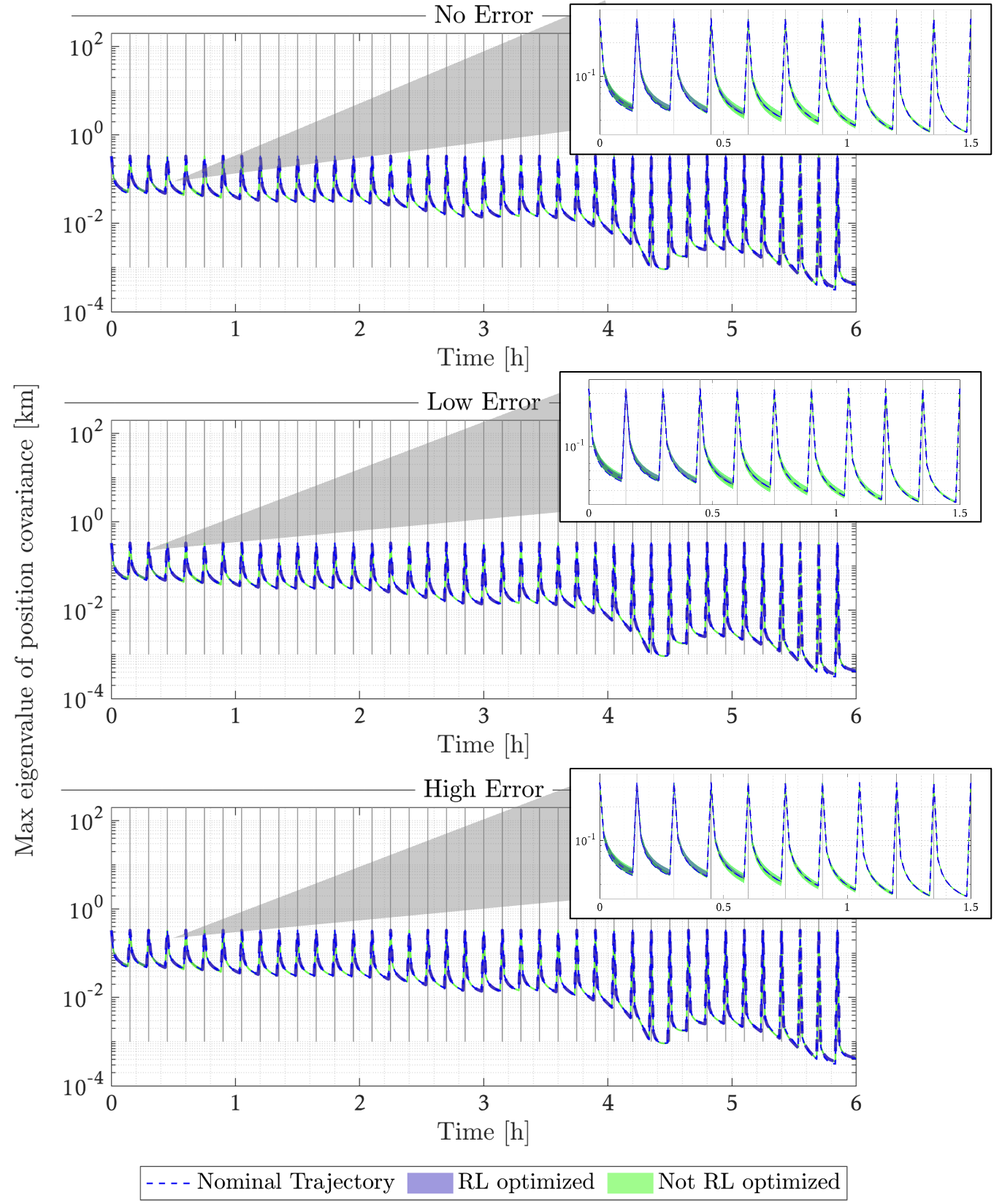}
    \caption{Observability: Distributions of the maximum position covariance over time under different error cases.}
    \label{eigcovlong}
\end{figure}

Fig. \ref{eigcovlong} shows the variation in the maximum position covariance over time for the MC trajectories under different error conditions, using an EKF. As the RL-trained result reaches the nominal trajectory around 2 h after $t_0$, significant deviations from the nominal are only observed within the maximum position eigenvalues at the start. The zoomed-in plots reveal that the RL-trained results maintain position covariances much closer to the nominal trajectory than other $\alpha$ strategies. This is especially evident for low and high $\Delta \overline{\bm{v}}$ error cases. This underscores the superior observability of the $\alpha_{RL}$ trajectories, even in the presence of thrust errors.

\begin{figure}[hbt!]
    \centering
    \includegraphics[width=0.8\linewidth]{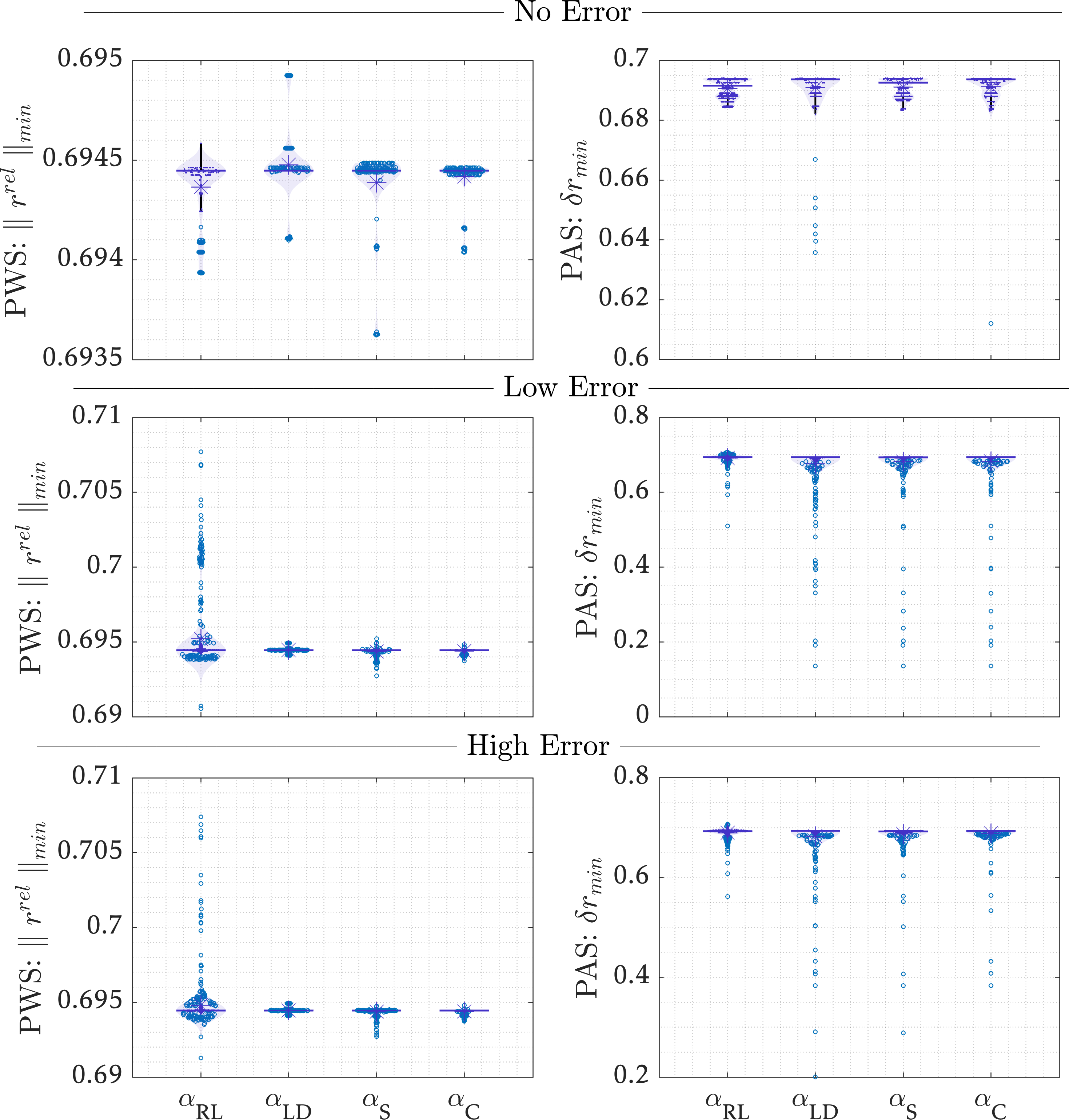}
    \caption{Safety: Left: Distribution of minimum distance to the target (for PWS safety). Right: Distribution of the minimum distance to the target when impulses are missed (for PAS safety).}
    \label{safelong}
\end{figure}

Fig. \ref{safelong} (left) shows the minimum distance between the target and the target in the MC simulations. This illustrates that the chaser does not breach the KOZ during its trajectory in any of the cases. However, Fig. \ref{safelong} (right) illustrates the minimum separation encountered when PAS safety is considered,  which shows that only the RL-trained controller does not break the KOZ barrier. Hence only  $\alpha_{RL}$ is both PWS and PAS safe. 

\begin{figure}[hbt!]
    \centering
    \includegraphics[width=0.8\linewidth]{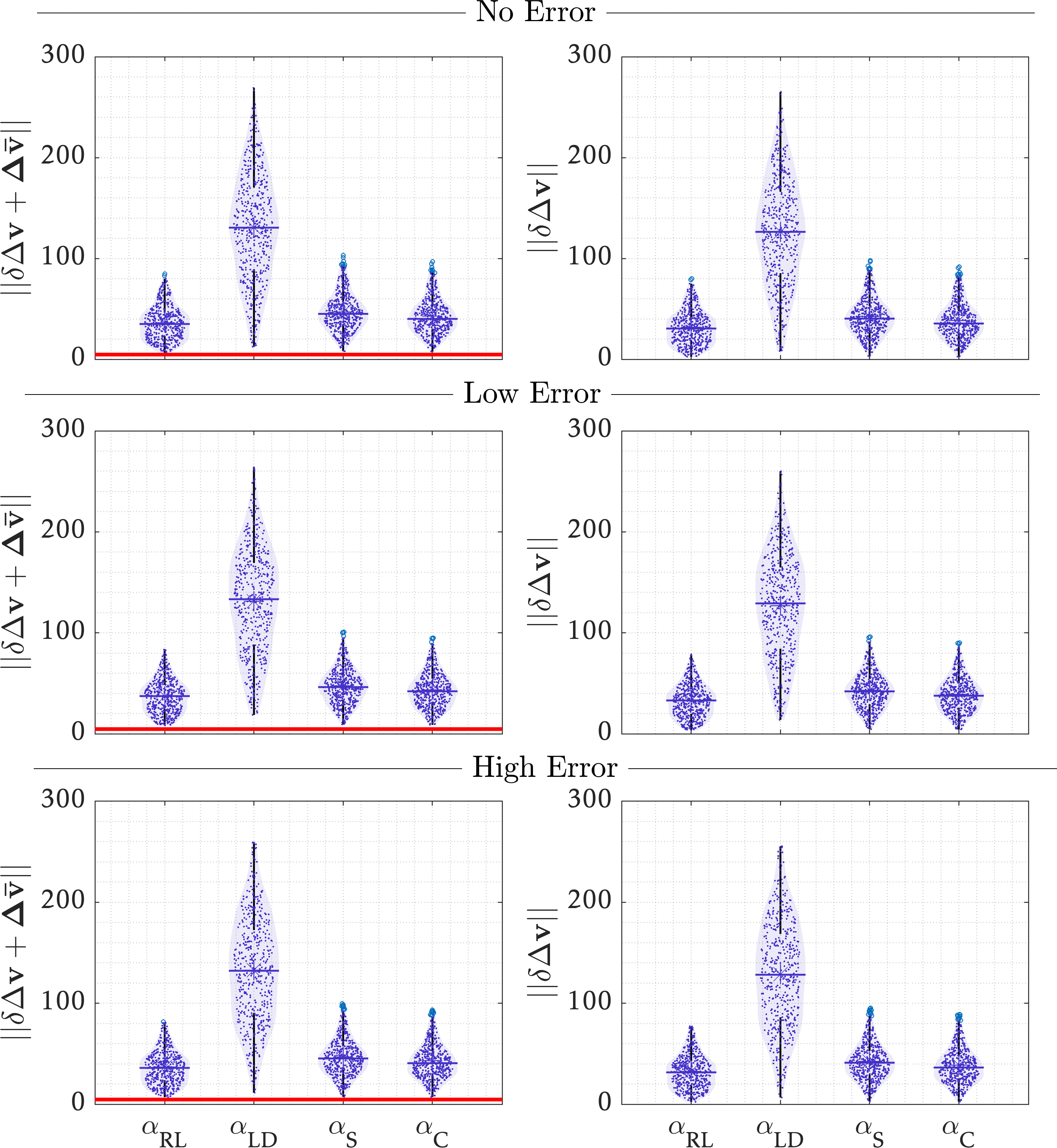}
    \caption{Comparison of the $\Delta v$ consumption of the RL optimized $\alpha$ with $\alpha$ derived by other means. The bold red line shows the total nominal $\Delta v$, denoted $\Delta \bar{v}$.}
    \label{dvlong}
\end{figure}

As seen in Fig. \ref{dvlong} and Table \ref{resultslong}, the performance of the RL guidance stands out in terms of its $\Delta v$ consumption.  According to Table \ref{resultslong}, $\alpha_{RL}$ consumes $16.58\%$ less $\Delta v$ compared to the next-best strategy, $\alpha_S$ in the case without errors. For the highest error case, this difference decreases slightly to $15.47\%$, but the RL results provide the lowest $\Delta v$ consumptions compared to other distributions even in the presence of high errors.  $\alpha_{RL}$ also has the best 99th percentile  ($P_{99}$) values for all error cases in this simulation.
The $\Delta v$ consumption of $\alpha_{RL}$ stands out in this manner likely because the RL reward function directly optimizes $\Delta v$ while only penalizing observability and safety if they fall below that of the nominal trajectory. In \cite{YUAN2022107812}, the authors report a mean $\Delta v$ consumption of 86.85 m/s, which is 2.36 times higher than the $\Delta v$ consumption obtained in the highest error case of this work.

\begin{table}[hbt!]
    \centering
    \begin{tabular}{ lll }
        \hline
       {$\alpha$} & \multicolumn{2}{c }{$\parallel \delta \Delta \pmb{v} + \Delta \pmb{v}\parallel$ {[}m/s{]}} \\ \cline{2-3}
                      & $\mu \pm \sigma$        & $[P_{25} , P_{75}, P_{99}]$   \\ \hline
                      & \multicolumn{2}{c }{No Error}                                \\ 
$\alpha_{RL}$             & $36.3814 \pm 17.0083$       & $[22.9455,46.4399,77.4094]$            \\
$\alpha_{LD}$             & $130.6799 \pm 53.6418$       & $[89.3742,170.1972,248.1801]$              \\
$\alpha_C$                & $42.9590 \pm 17.9896$       & $[30.0029,52.3218,88.4102]$           \\
$\alpha_S$                & $46.9150 \pm 18.8545$       & $[34.0597,57.1734,94.144]$               \\
\hline

                      & \multicolumn{2}{c }{Low Error}  \\ 
$\alpha_{RL}$             & $37.4193 \pm 16.5074$       & $[24.3864,48.5012,78.0634]$              \\
$\alpha_{LD}$             & $131.1536 \pm 52.3019$       & $[88.3215,169.1716,249.5346]$               \\
$\alpha_C$                & $43.6415 \pm 17.5211$       & $[30.3208,54.3998,89.0115]$               \\
$\alpha_S$                & $47.5103 \pm 18.3366$       & $[33.8143,58.8386,94.8156]$                  \\
\hline           & \multicolumn{2}{c }{High Error} \\ 
$\alpha_{RL}$             & $36.8001 \pm 16.8291$       & $[23.7529,46.795,80.2011]$            \\
$\alpha_{LD}$             & $132.7175 \pm 52.8530$       & $[90.4364,172.3897,253.9687]$              \\
$\alpha_C$                & $42.9698 \pm 17.8540$       & $[30.0495,52.728,91.4231]$               \\
$\alpha_S$                & $46.8578 \pm 18.6912$       & $[33.8758,57.2557,97.3885]$              \\
\hline
        \end{tabular}

    \caption{Comparison of the $\Delta v$ consumption using Monte-Carlo simulations under different error levels}
    \label{resultslong}
\end{table}

 \section{Conclusion}

This work develops an approach for the design and guidance of far-range operations under angles-only navigation, considering safety, robustness, and observability. The developed approach consists of two stages. Firstly, particle swarm optimization is used to develop a nominal, multi-impulse far-range trajectory that is both PWS and PAS safe and observable. Observability is ensured by adjusting the impulses imparted to optimize the angles-only measurements and using nonlinear Keplerian dynamics in the problem formulation. Point-wise safety is assured by maintaining a minimum separation between the target and the chaser throughout the transfer. PAS safety is implemented by maintaining a minimum distance between the chaser and the target for one target orbital period if an impulse is missed. Secondly, an RL-based guidance scheme is developed for spacecraft guidance, which aims to minimize fuel consumption while maintaining observability and safety. The constraint satisfaction challenges of particle swarm optimization and RL are alleviated through the problem formulation, which incorporates the Lambert method to guarantee that the target state is always reached. The guidance controller is trained on a distribution of initial states and its performance is evaluated via a 500-sample Monte Carlo simulation with initial state variations and nominal $\Delta \overline{v}$ errors. The observability of the MC results is evaluated using an Extended Kalman filter, and safety is investigated by studying the PWS and PAS minimum separation distances observed by the Monte Carlo samples. The results show that the RL controller consumes significantly less fuel
than other benchmark distributions while complying with safety and observability requirements. In the test case studied, the RL controller is seen to overcome the challenges of angles only navigation to generate observable and safe trajectories with the aid of the particle swarm-optimized nominal trajectory.

\section*{Acknowledgments}
The author(s) wish to acknowledge the use of New Zealand eScience Infrastructure (NeSI) high performance computing facilities, consulting support and/or training services as part of this research. New Zealand's national facilities are provided by NeSI and funded jointly by NeSI's collaborator institutions and through the Ministry of Business, Innovation \& Employment's Research Infrastructure programme. URL \url{https://www.nesi.org.nz}.

\bibliographystyle{unsrt}

\end{document}